\documentclass[12pt]{amsart}
\usepackage{a4wide,amsfonts,amssymb,stmaryrd,amscd,amsmath,latexsym,amsbsy}
\numberwithin{equation}{section}
\theoremstyle{plain}
\newtheorem{thm}{Theorem}[section]

\newtheorem{prop}[thm]{Proposition}

\newtheorem{defi}[thm]{Definition}
\newtheorem{lem}[thm]{Lemma}
\newtheorem{cor}[thm]{Corollary}

\theoremstyle{remark}
\newtheorem{rema}[thm]{Remark}
\newtheorem{eg}[thm]{{Example}}
\title{Periodic integrable systems with delta-potentials}
\author{E. Emsiz, E.M. Opdam and J.V. Stokman}
\address{KdV Institute for Mathematics, Universiteit van Amsterdam,
Plantage Muidergracht 24, 1018 TV Amsterdam, The Netherlands.}
\email{eemsiz@science.uva.nl,
opdam@science.uva.nl,
jstokman@science.uva.nl}
\begin{document}
\begin{abstract}
In this paper we study root system generalizations of
the quantum Bose-gas on the circle with pair-wise delta-function interactions.
The underlying symmetry structures are shown to be governed by
the associated graded algebra of Cherednik's (suitably filtered)
degenerate double affine Hecke algebra,
acting by Dunkl-type differential-reflection
operators. We use Gutkin's generalization of the equivalence between
the impenetrable Bose-gas and the free Fermi-gas to
derive the Bethe ansatz equations and the Bethe ansatz eigenfunctions.
\end{abstract}
\maketitle

\section{Introduction}\label{sect1}
Given any affine root system $\Sigma$, Gutkin and Sutherland
\cite{GS}, \cite{Su} defined a quantum integrable system whose Hamiltonian $-\Delta+\mathcal{V}$ has a potential
$\mathcal{V}$ expressible as a weighted sum of delta-functions at the affine root hyperplanes of $\Sigma$.
For the affine root system of type $A$, the quantum integrable system essentially reduces to the
quantum Bose-gas on the circle with pair-wise delta-function interactions,
which has been subject of intensive studies over the past 40 years.

The special case of the impenetrable Bose-gas on the circle
was exactly solved by relating the model to the free Fermi-gas on the circle (see Girardeau \cite{Gi}).
Soon afterwards fundamental progress was made for arbitrary pair-wise delta
function interactions by Lieb \& Liniger \cite{LL}, Yang \cite{Y} and
Yang \& Yang \cite{YY}, leading to the derivation of
the associated Bethe ansatz equations and Bethe ansatz eigenfunctions. Yang \& Yang
\cite{YY} showed that the solutions of the Bethe ansatz equations
are controlled by a strictly convex master function.
One of the aims of the present paper is to generalize these
results to Gutkin's and Sutherland's quantum integrable systems associated to affine root
systems.

Quantum Calogero-Moser systems
are root system generalizations of quantum Bose-gases on the line or circle
with long range pair-wise interactions. In special cases
quantum Calogero-Moser systems naturally arose from harmonic analysis on symmetric spaces.
A decisive role in the studies of quantum Calogero-Moser systems has been played by certain
non-bosonic analogs of these systems, which are defined in terms of
Dunkl-type commuting differential-reflection operators. Suitable degenerations of affine Hecke algebras
naturally appear here as the fundamental objects
governing the algebraic relations between the Dunkl-type operators and the natural Weyl group action.

In this paper we define Dunkl-type commuting differential-reflection operators
associated to the root system generalizations of the quantum Bose-gas with
delta-function interactions. We furthermore show that the Dunkl-type
operators, together with the natural affine Weyl group action,
realize a faithful representation of the
associated graded algebra of Cherednik's \cite{C} (suitably filtered)
degenerate double affine Hecke algebra. These results show that these quantum integrable systems
naturally fit into the class of quantum Calogero-Moser integrable systems, a point of view
which also has been advertised from the perspective of harmonic analysis in \cite[Sect. 5]{HO}.

The quantum integrable systems under consideration for affine root systems $\Sigma$ of classical type
still have reasonable physical interpretations in terms of interacting one-dimensional quantum bosons.
In these cases various results of the present paper can be found in the vast physics literature on this
subject. We will give the precise connections to the literature in the main body of the text.

The knowledge on the quantum Bose-gas with pair-wise delta-function interactions
still far exceeds the knowledge on its root system generalizations.
In fact, an important feature of the quantum Bose-gas with pair-wise delta-function interactions
is its realization as the restriction to a fixed particle sector of the quantum integrable
field theory in $1+1$ dimensions governed by the quantum nonlinear
Schr{\"o}dinger equation.
This point of view has led to the study of this model
by quantum inverse scattering methods.
With these methods a proof of full orthogonality of the Bethe
eigenfunctions on a period box (with respect to Lebesgue measure)
is derived in \cite{D} and the quadratic norms of the Bethe
eigenfunctions are evaluated in terms of the determinant of the Hessian of the master function
(conjectured by Gaudin \cite[Sect. 4.3.3]{Ga} and proved by Korepin \cite{K}).

At this point we can only speculate on the generalizations of these results
to arbitrary root systems. The quantum inverse scattering
techniques are only in reach for classical root systems,
in which case we have quantum field theories with (non)periodic
integrable boundary conditions to our disposal, see \cite{Sk}.
In general it seems reasonable to expect that
the Bethe eigenfunctions are orthogonal on a fundamental domain for the reflection representation
of the affine Weyl group (with respect to Lebesgue measure),
and that their quadratic norms are expressible in terms of the determinant of the Hessian of the master function
at the associated spectral point.

The contents of the paper is as follows. Sections \ref{sect2} and \ref{sect3}
are meant to introduce the quantum integrable systems and to state and clarify the
results on the associated spectral problem. We first introduce in Section \ref{sect2}
the relevant notations on affine root systems. Following Gutkin \cite{G} we formulate the spectral problem for the
quantum integrable systems under consideration as an explicit boundary value problem.
We state the main results on the boundary value problem (Bethe ansatz equations and Bethe ansatz eigenfunctions)
and we introduce the associated master function. In Section \ref{sect3} we formulate the
analog of Girardeau's
equivalence between the impenetrable Bose-gas and the free Fermi-gas on the circle for the quantum integrable systems
under consideration.

In Section \ref{sect4} we introduce Dunkl-type commuting differential-reflection operators and
show that they realize, together with the natural affine Weyl group action, a faithful realization of the associated
graded algebra $H$ of Cherednik's \cite{C} (suitably filtered) degenerate double affine Hecke algebra. In Section \ref{sect5} we show that
Gutkin's \cite{G} integral-reflection operators, together with the ordinary directional derivatives, yield an
equivalent realization of $H$. The equivalence is realized by Gutkin's \cite{G} propagation operator.
We furthermore show that the Dunkl operators naturally act on a space of functions with higher order
normal derivative jumps over the affine root hyperplanes.

In Section \ref{sect6} we return to the boundary value problem of Section \ref{sect2}. Using the Hecke-type algebra $H$
we refine and clarify Gutkin's \cite{G} generalization of Girardeau's equivalence between the boundary value problem for
the impenetrable Bose gas and the boundary value problem for the free Fermi-gas as formulated in Section \ref{sect3}.
The results in this section entail that
the boundary value problem is equivalent to a boundary value problem with trivial boundary value conditions,
at the cost of having to deal with a non-standard affine Weyl group action.
In Section \ref{sect7} we study the reformulated boundary value problem, leading in Section \ref{sect8}
to the derivation of the Bethe ansatz equations. In Section \ref{sect9} we study the master function and show
how it leads to a natural parametrization of the solutions of the Bethe ansatz equations. In Section \ref{sect10}
the solutions of the Bethe ansatz equations are further analyzed. In Section \ref{sect11} it is proved
that the boundary value problem has solutions if and only if the associated spectral value is a {\it regular}
solution of the Bethe ansatz equations. In case of root system of type $A$, this is known as the Pauli principle
for the interacting bosons.
\vspace{.2cm}\\
{\bf Acknowledgments:}
This research was supported in part by a Pionier grant
of the Netherlands organization for Scientific Research (NWO)
(Emsiz and Opdam (program coordinator)).
Stokman was supported by the Royal Netherlands
Academy of Arts and Sciences (KNAW) and by the Netherlands Organization
for Scientific Research (NWO) in the VIDI-project ``Symmetry and
modularity in exactly solvable models''.


\section{The boundary value problem}\label{sect2}
In this section we recall Gutkin's \cite{G} reformulation of the
spectral problem for periodic integrable systems with delta-potentials
in terms of a concrete boundary value problem. We furthermore state the
main results on the solutions of the boundary value problem and we detail the physical background.

In order to fix notations we start by recalling some well known facts on
affine root systems, see e.g. \cite{Hu} for a detailed exposition.
Let $V$ be an Euclidean space of dimension $n$. Let $\Sigma_0$ be a finite, irreducible crystallographic
root system in the dual Euclidean space $V^*$. We denote $\langle \cdot,\cdot\rangle$ for the inner product on $V^*$
and $\|\cdot\|$ for the corresponding norm.
The co-root of $\alpha\in \Sigma_0$ is the unique vector $\alpha^\vee\in V$ satisfying
\[\xi(\alpha^\vee)=
\frac{2\langle\xi,\alpha\rangle}{\|\alpha\|^2},\quad \forall\,\xi\in V^*.
\]
We write $\Sigma_0^\vee=\{\alpha^\vee\}_{\alpha\in\Sigma_0}$ for the resulting co-root system in $V$.
We fix a basis $I_0=\{a_1,\ldots,a_n\}$ for the root system $\Sigma_0$.
Let $\Sigma_0=\Sigma_0^+\cup \Sigma_0^-$ be the
corresponding decomposition in positive and negative roots. We
denote $\rho\in V^*$ for the half sum of positive roots and $\varphi\in \Sigma_0^+$
for the highest root with respect to the basis $I_0$. The highest root $\varphi$ is a long root
in $\Sigma_0$. We define the fundamental Weyl
chamber in $V^*$ by
\begin{equation}\label{chamber}
V_+^*=\{\xi\in V^* \, | \, \xi(\alpha^\vee)>0\quad
\forall\,\,\alpha\in\Sigma_0^+ \}.
\end{equation}

Let $\widehat{V}$ be the vector space of affine linear functionals on $V$.
Then $\widehat{V}\simeq V^*\oplus \mathbb{R}$ as vector spaces,
where the second component
is identified with the constant functions on $V$.
The gradient map $D: \widehat{V}\rightarrow V^*$
is the projection onto $V^*$ along this decomposition.

The subset $\Sigma=\Sigma_0+\mathbb{Z}\subset \widehat{V}$
is the affine root system associated to $\Sigma_0$. We extend
the basis $I_0$ of $\Sigma_0$ to a basis
$I=\{a_0=-\varphi+1,a_1,\ldots,a_n\}$ of the affine root system
$\Sigma$. Observe that $D$ maps $\Sigma$ onto $\Sigma_0$.

For a root $a\in \Sigma$,
\[s_a(v)=v-a(v)Da^\vee,\qquad v\in V
\]
defines the orthogonal reflection in the root hyperplane $V_a:=a^{-1}(0)$.
The affine Weyl group $W$ associated to $\Sigma$ is the sub-group of the affine linear isomorphisms of $V$
generated by the orthogonal reflections $s_a$ ($a\in \Sigma$). The sub-group $W_0\subset W$
generated by the orthogonal reflections $s_\alpha$ ($\alpha\in \Sigma_0$)
is the Weyl group associated to $\Sigma_0$. We denote $w_0$
for the longest Weyl group element in $W_0$.
It is well known that $W$ (respectively $W_0$) is a Coxeter group with Coxeter generators
the simple reflections $s_j=s_{a_j}$ for $j=0,\ldots,n$
(respectively $s_j$ for $j=1,\ldots,n$).

A second important presentation of $W$ is given by
\begin{equation}\label{secondformW}
W\simeq W_0\ltimes Q^\vee,
\end{equation}
with $Q^\vee=\mathbb{Z}\Sigma_0^\vee\subset V$ the co-root lattice of
$\Sigma_0$, acting by translations on $V$. The gradient map $D$ induces
a surjective group homomorphism $D: W\rightarrow W_0$ by $D(s_a)=s_{Da}$
for $a\in\Sigma$. Alternatively, $Dw=v$ if $v\in W_0$ is the
$W_0$-component of $w$ in the semi-direct product decomposition
\eqref{secondformW}.

The space $\widehat{V}$ of affine linear functionals on $V$ is a $W$-module by
$(wf)(v)=f(w^{-1}v)$ ($w\in W, f\in \widehat{V}, v\in V$). Observe that
$V^*$ is $W_0$-stable, and
\[s_\alpha(\xi)=\xi-\xi(\alpha^\vee)\alpha,\qquad \xi\in V^*
\]
for roots $\alpha\in \Sigma_0$. Furthermore,
\[s_\alpha(\Sigma_0)=\Sigma_0,\qquad s_a(\Sigma)=\Sigma
\]
for $\alpha\in\Sigma_0$ and $a\in\Sigma$. The length of $w\in W$ is defined by
$l(w)=\#\bigl(\Sigma^+\cap w^{-1}\Sigma^-\bigr)$. Alternatively, $l(w)$ is the
minimal positive integer $r$ such that $w\in W$ can be written as
product of $r$ simple reflections. Such an expression $w=s_{j_1}s_{j_2}\cdots s_{j_{l(w)}}$
($j_k\in\{0,\ldots,n\}$) is called reduced.

The weight lattice of $\Sigma_0$ is defined by
\[
P=\{\lambda\in V^* \,\,\, | \,\,\,
\lambda(\alpha^\vee)\in\mathbb{Z}\quad \forall\, \alpha\in\Sigma_0 \}.
\]
Another convenient description is
\begin{equation}\label{Palt}
P=\{\lambda\in V^* \,\,\, | \,\,\,
w\lambda(\varphi^\vee)\in\mathbb{Z} \quad \forall\,w\in W_0 \},
\end{equation}
which follows from the fact that $Q^\vee$ is already spanned
over $\mathbb{Z}$ by the short co-roots in $\Sigma_0^\vee$.
We denote $P^+$ (respectively $P^{++}$) for the cone of
dominant (respectively strictly dominant) weights with respect to
the choice $\Sigma_0^+$ of positive roots in $\Sigma_0$. Recall
that $P^{++}=\rho+P^+$.

We write $V_{irreg}=\bigcup_{a\in\Sigma^+}V_a$ for the irregular vectors in $V$ with respect
to the affine root hyperplane arrangement $\{V_a\,\, | \,\,  a\in\Sigma^+\}$. Its open, dense complement $V_{reg}:=V\setminus V_{irreg}$ is called
the set of regular vectors in $V$.

We denote $\mathcal{C}$ for the collection of connected components of $V_{reg}$. An element $C\in \mathcal{C}$
is called an alcove. The affine Weyl group $W$ acts simply transitively on $\mathcal{C}$.
Explicitly, $V_{reg}=\bigcup_{w\in W}w(C_+)$ (disjoint union)
with the fundamental alcove $C_+$ defined by
\[C_+=\{\,v\in V \, | \, a_j(v)>0\,\,\, (j=0,\ldots,n)\,\}.
\]
We call a vector $v\in V_a$ ($a\in\Sigma^+$) sub-regular if it does not
lie on any other root hyperplane $V_b$ ($a\not=b\in\Sigma^+$).

The symmetric algebra $S(V)$ is canonically a $W_0$-module
algebra. Using the standard identification $S(V)\simeq P(V^*)$
where $P(V^*)$ is the algebra of real-valued polynomial functions on $V^*$,
the $W_0$-module structure takes the form
\[(wp)(\xi)=p(w^{-1}\xi),\qquad w\in W_0,\,\,\xi\in V^*.
\]
We denote $S(V)^{W_0}$ and $P(V^*)^{W_0}$
for the subalgebra of $W_0$-invariants in $S(V)$ and $P(V^*)$,
respectively.

Let $\partial_v$ ($v\in V$) be the derivative in direction $v$,
\[(\partial_vf)(u)=\frac{d}{dt}\bigg|_{t=0}f(u+tv)
\]
for $f$ continuously differentiable at $u\in V$. The assignment $v\mapsto \partial_v$
uniquely extends to an algebra isomorphism of $S(V)$ onto the
algebra of constant coefficient differential operators on $V$ (say acting on $C^\infty(V)$). We
denote $p(\partial)$ for the constant coefficient differential
operator corresponding to $p\in S(V)\simeq P(V^*)$. For example,
the $W_0$-invariant constant coefficient differential operator $p_2(\partial)$
associated to the polynomial $p_2(\cdot)=\|\cdot\|^2\in P(V^*)^{W_0}$
is the Laplacean $\Delta$ on $V$.

The quantum integrable system which we will define now in a moment
depends on certain coupling constants
called multiplicity functions.

\begin{defi}
A multiplicity function $k$ is a $W$-invariant function $k: \Sigma\rightarrow \mathbb{R}$
satisfying $k(a)=k(Da)$ for all $a\in \Sigma$.
\end{defi}

Unless stated explicitly otherwise, we fix a strictly positive multiplicity function
$k:\Sigma \rightarrow \mathbb{R}_{>0}$. To simplify notations
we write $k_a$ for the value of $k$ at the root $a\in \Sigma$.

We define the quantum Hamiltonian $\mathcal{H}_k$ by
\begin{equation}\label{Hamiltonian}
\mathcal{H}_k=-\Delta+\sum_{a\in\Sigma}k_a\delta(a(\cdot)),
\end{equation}
where $\delta$ is the Kronecker delta-function.
Here we interpret $\mathcal{H}_k$ as a linear map
$\mathcal{H}_k: C(V)\rightarrow D^\prime(V)$,
with $C(V)$ the complex-valued continuous functions on $V$ and
$D^\prime(V)$ the space of distributions on $V$, as
\begin{equation}\label{distribution}
\bigl(\mathcal{H}_kf\bigr)(\phi):=-\int_Vf(v)\bigl(\Delta\phi\bigr)(v)dv
+\sum_{a\in\Sigma}\frac{k_a}{\|Da^\vee\|}\int_{V_a}f(v)\phi(v)d_av
\end{equation}
for a test function $\phi$, with $dv$ the Euclidean volume measure on $V$
and $d_av$ ($a\in\Sigma^+$) the corresponding volume measure on the root
hyperplane $V_a$.

The quantum Hamiltonian $\mathcal{H}_k$ and the
associated quantum physical system has been studied in e.g.
\cite{Su}, \cite{GS} and \cite{G}. A key step in these
investigations is the reformulation of the spectral problem for $\mathcal{H}_k$
in terms of an explicit boundary value problem for the Laplacean $\Delta$ on $V$,
which we now proceed to recall.

Let $CB^1(V)$ be the space of complex valued continuous functions $f$ on $V$
whose restriction $f|_C$ to an alcove $C\in\mathcal{C}$ has a
continuously differentiable
extension to some open neighborhood $\widetilde{C}\supset
\overline{C}$.
Let $C^{1,(k)}(V)$ be the space of functions $f\in CB^1(V)$ which satisfy the derivative jump conditions
\begin{equation}\label{jumpcond}
\bigl(\partial_{Da^\vee}f\bigr)(v+0Da^\vee)-\bigl(\partial_{Da^\vee}f\bigr)(v-0Da^\vee)=
2k_af(v)
\end{equation}
for sub-regular vectors $v\in V_a$ \textup{(}$a\in\Sigma^+$\textup{)}.

\begin{prop}\label{visual1}
For $f\in CB^1(V)$ and $E\in\mathbb{C}$ the following two statements are equivalent.\\
{\bf (i)} $\mathcal{H}_kf=Ef$ as distributions on $V$.\\
{\bf (ii)} $f\in C^{1,(k)}(V)$ and $\Delta f|_{V_{reg}}=-Ef|_{V_{reg}}$
as distributions on $V_{reg}$.\\
A function $f\in CB^1(V)$ satisfying these equivalent conditions
is smooth on $V_{reg}$.
\end{prop}
\begin{proof}
The first part of the proposition follows from a straightforward application
of Green's identity (cf. the proof of \cite[Thm. 2.7]{G}).
The last statement follows from the fact that
the constant coefficient differential operator $\Delta+E$ on $V$ is (hypo)elliptic.
\end{proof}
The quantum physical system with quantum Hamiltonian $\mathcal{H}_k$ is
known to be integrable. The common spectral problem for the
associated quantum conserved integrals has been translated by Sutherland
and Gutkin \cite{Su}, \cite{GS} into the following boundary value problem.
\begin{defi}\label{BVP}
Fix a spectral parameter $\lambda\in V_{\mathbb{C}}^*:=\mathbb{C}\otimes_{\mathbb{R}}V^*$.
We denote $\textup{BVP}_k(\lambda)$ for the space of functions $f\in C^{1,(k)}(V)$
solving \textup{(}in distributional sense\textup{)} the system
\begin{equation}\label{diffeqn}
p(\partial)f\big|_{V_{reg}}=p(\lambda)f\big|_{V_{reg}}\qquad \forall \, p\in
S(V)^{W_0}
\end{equation}
of differential equations away from the root hyperplane configuration $\bigcup_{a\in\Sigma^+}V_a$.
\end{defi}
\begin{rema}\label{ellarg}
Since $\Delta=p_2(\partial)$ is the Laplacean on $V$, Proposition \ref{visual1} implies that a
function $f\in \hbox{BVP}_k(\lambda)$ is smooth on $V_{reg}$ and satisfies the differential
equations \eqref{diffeqn} in the strong sense. The fact that $f$ is an eigenfunction of all
$W_0$-invariant constant coefficient differential operators on $V_{reg}$
in fact implies that $f|_C$ is the restriction of a (necessarily unique) analytic function
on $V$ for all alcoves $C\in\mathcal{C}$, see \cite{S}.
\end{rema}

The central theme of this paper is the study of the subspace
$\textup{BVP}_k(\lambda)^W\subset \textup{BVP}(\lambda)$
of $W=W_0\ltimes Q^\vee$-invariant solutions, where $W$ acts on
$\textup{BVP}_k(\lambda)\subset C^{1,(k)}(V)$ by
\begin{equation}\label{usual}
(wf)(v)=f(w^{-1}v)
\end{equation}
for $w\in W$ and $v\in V$.
Our focus on $W$-invariant
solutions thus amounts to studying the bosonic (=$W_0$-invariant) theory
of the quantum system under $Q^\vee$-periodicity constraints
(or equivalently, we view the quantum system on the torus $V/Q^\vee$).
\begin{eg}[Free case $k\equiv 0$]\label{k=nul}
A function $f\in \textup{BVP}_0(\lambda)$ is a distribution solution of the (hypo)elliptic
constant coefficient differential operator $\Delta-p_2(\lambda)$ on $V$,
hence $f$ is smooth on $V$ (cf. Proposition \ref{visual1}). Combined with Remark \ref{ellarg}
we conclude that a function $f\in \textup{BVP}_0(\lambda)$ is analytic on $V$.
Then $\textup{BVP}_0(\lambda)^W$ ($\lambda\in V_{\mathbb{C}}^*$) are the common
eigenspaces of the quantum conserved integrals for
the free bosonic quantum integrable system on $V/Q^\vee$
associated to the Laplacean $\Delta$ on $V$. It is easy to show that
$\textup{BVP}_0(\lambda)^W$ is zero-dimensional unless $\lambda\in 2\pi i P$,
in which case it is spanned by the plane wave
\[\phi_\lambda^0=\frac{1}{\#W_0}\sum_{w\in W_0}e^{w\lambda}
\]
\textup{(}cf. the analysis in the
impenetrable case $k\equiv\infty$ in Section \ref{sect3}\textup{)}.
\end{eg}

The quantum Hamiltonian \eqref{Hamiltonian} for $\Sigma_0$ of
type $A_n$ takes the explicit form
\[-\Delta+k\sum_{m\in\mathbb{Z}}\,\sum_{1\leq i\not=j\leq
n+1}\delta(x_i-x_j+m).
\]
Here we have embedded $V$ into $\mathbb{R}^{n+1}$ as the hyperplane defined by
$x_1+\cdots +x_{n+1}=0$. The study of $W$-invariant
solutions to the boundary value problem then
essentially amounts to analyzing the
spectral problem for the system describing $n+1$ quantum bosons
on the circle with
pair-wise repulsive delta-function interactions.
In this special case the quantum system
has been extensively studied in the physics literature,
see e.g. \cite{Gi}, \cite{LL}, \cite{Y}, \cite{YY}, \cite{Ga}
and \cite{IK}. The upgrade to other classical root systems
amounts to adding particular reflection terms to the physical model,
see e.g. \cite{Sk},
\cite{Ch0}, \cite{Ga}, \cite{H}, \cite{KH} and \cite{MW}.

We are now in a position to formulate the main results on
the solution space of the boundary value
problem. We call the spectral value
$\lambda\in V_{\mathbb{C}}^*=V^*\oplus iV^*$
regular if its isotropy sub-group in $W_0$ is trivial
(equivalently, $\lambda(\alpha^\vee)\not=0$ for all $\alpha\in\Sigma_0$).
We call $\lambda$ singular
otherwise. Furthermore, $\lambda$ is called real
(respectively purely imaginary)
if $\lambda\in V^*$ (respectively $\lambda\in iV^*$).
Define the $c$-function by
\begin{equation}\label{cfunction}
\widetilde{c}_k(\lambda)=
\prod_{\alpha\in\Sigma_0^+}
\frac{\lambda(\alpha^\vee)+k_\alpha}{\lambda(\alpha^\vee)}
\end{equation}
as rational function of $\lambda\in V_{\mathbb{C}}^*$,
cf. \cite{Ga}, \cite{HO}.
\begin{thm}\label{main}
Let $\lambda\in V_{\mathbb{C}}^*$. The space $\textup{BVP}_k(\lambda)^W$
of $W$-invariant solutions to the boundary value problem is
one-dimensional or zero-dimensional.
It is one-dimensional if and only if the spectral value $\lambda$ is
a purely imaginary,
regular solution of the Bethe ansatz equations
\begin{equation}\label{BAE}
e^{w\lambda(\varphi^\vee)}=\prod_{\alpha\in\Sigma_0^+}
\left(\frac{w\lambda(\alpha^\vee)-k_\alpha}{w\lambda(\alpha^\vee)+k_\alpha}
\right)^{\alpha(\varphi^\vee)}\qquad \forall\, w\in W_0.
\end{equation}
If $\textup{BVP}_k(\lambda)^W$ is one-dimensional, then there
exists a unique $\phi_\lambda^k\in \textup{BVP}_k(\lambda)^W$
normalized by $\phi_\lambda^k(0)=1$. The solution $\phi_\lambda^k$
is the unique $W$-invariant function satisfying
\begin{equation}\label{BAF}
\phi_\lambda^k(v)=\frac{1}{\#W_0}\sum_{w\in
W_0}\widetilde{c}_k(w\lambda)e^{w\lambda(v)},\qquad v\in
\overline{C_+}.
\end{equation}
\end{thm}
We give a reformulation of Theorem \ref{main} in Section \ref{sect3}.
The Bethe ansatz equations are derived in Section \ref{sect8}.
The regularity constraint on $\lambda$ is proved in Section \ref{sect11}.
\begin{rema}\label{BAEremark}
The Bethe ansatz equations \eqref{BAE} can be rewritten as
\begin{equation}\label{BAEconventional}
e^{w\lambda(\varphi^\vee)}=\frac{w\lambda(\varphi^\vee)-k_\varphi}
{w\lambda(\varphi^\vee)+k_\varphi}\prod_{\alpha\in\Sigma_0^+\cap
s_{\varphi}\Sigma_0^-}\frac{w\lambda(\alpha^\vee)-k_\alpha}
{w\lambda(\alpha^\vee)+k_\alpha},
\qquad \forall\, w\in W_0,
\end{equation}
due to the fact that for $\alpha\in \Sigma_0^+$,
\begin{equation}\label{highestrootproperty}
\alpha(\varphi^\vee)=
\begin{cases}
2\qquad &\hbox{ if } \alpha=\varphi,\\
1\qquad &\hbox{ if } \alpha\in\bigl(\Sigma_0^+\cap
s_\varphi\Sigma_0^-\bigr)\setminus \{\varphi\},\\
0\qquad &\hbox{ if } \alpha\in\Sigma_0^+\cap s_\varphi\Sigma_0^+.
\end{cases}
\end{equation}
\end{rema}

A key role in the analysis of the Bethe ansatz equations \eqref{BAE} is played
by the following {\it master function}.
\begin{defi}
The master function $S_k:P\times V^*\rightarrow \mathbb{R}$ is
defined by
\begin{equation}\label{master}
S_k(\mu,\xi)=\frac{1}{2}\|\xi\|^2-2\pi\langle \mu,\xi\rangle
+\frac{1}{2}\sum_{\alpha\in \Sigma_0}\|\alpha\|^2\int_0^{\xi(\alpha^\vee)}
\arctan\left(\frac{t}{k_\alpha}\right)dt.
\end{equation}
\end{defi}
The master function $S_k$ enters into the description
of the set $\hbox{BAE}_k$ of solutions $\lambda\in iV^*$
of the Bethe ansatz equations \eqref{BAE} in the following way.
\begin{prop}\label{parametrization}
For $\mu\in P$ there exists a unique extremum $\widehat{\mu}_k\in V^*$ of
the master function $S_k(\mu,\cdot)$.
The assignment $\mu\mapsto i\widehat{\mu}_k$ defines a $W_0$-equivariant
bijection
$P\overset{\sim}{\longrightarrow} \textup{BAE}_k$.
\end{prop}
The proof of Proposition \ref{parametrization},
which hinges on the strict convexity of $S_k(\mu,\cdot)$
($\mu\in P$), is given in Section \ref{sect9}.
The regularity condition on the spectrum in Theorem \ref{main}
also turns out to be a consequence of
the strict convexity of the master function
$S_k(\mu,\cdot)$ ($\mu\in P$), see Section \ref{sect11}.

The following proposition yields precise information on the location of
the deformed weight $i\widehat{\mu}_k\in \hbox{BAE}_k$.
\begin{prop}\label{mg}
For $\mu\in P^{+}$ and $\beta\in\Sigma_0^+$ we have
\begin{equation}\label{momentgaps}
\frac{2\pi\mu(\beta^\vee)}{\left(1+\frac{h_k}{n}\right)}\leq
\widehat{\mu}_k(\beta^\vee)\leq 2\pi\mu(\beta^\vee),
\end{equation}
where $h_k=2\sum_{\alpha\in\Sigma_0}k_\alpha^{-1}$. Furthermore,
$\mu\in P^+$ if and only if $\widehat{\mu}_k\in \overline{V_+^*}$.
\end{prop}
Proposition \ref{mg} is proved in Section \ref{sect10}.
The lower bound in \eqref{momentgaps} shows how far away the spectral values
$\widehat{\mu}_k\in V_{+}^*$ ($\mu\in P^{++}$) are from being
singular.

The Bethe ansatz functions $\phi_\lambda^k$ and the necessity of the
Bethe ansatz equations \eqref{BAE} on the allowed spectrum were obtained by Lieb and Liniger
\cite{LL} for root system $\Sigma_0$ of type $A_n$, and soon after generalized to root system $\Sigma_0$ of type $D_n$
by Gaudin \cite{Ga0}, \cite{Ga} (see also \cite{KH}). For $\Sigma_0$ of type $A_n$, Yang and Yang \cite{YY} introduced the master
function $S$ (also known as the Yang-Yang action) and derived the special case of Proposition \ref{parametrization} using
its strict convexity.

In physics literature the regularity of the spectral parameter
$\lambda$ (see Theorem \ref{main}) is usually imposed as an additional requirement, since it
automatically ensures that eigenstates admit a plane wave expansion within any alcove $C\in \mathcal{C}$.
The regularity condition for root system $\Sigma_0$ of type $A_n$ can be viewed
as a Pauli type principle for the interacting quantum {\it bosons}, since it
implies that the momenta of the quantum bosons are pair-wise different.
An actual proof of the regularity of the spectrum
was obtained by Izergin and Korepin \cite{IK} using quantum inverse
scattering methods. In this derivation the regularity condition again follows from
the strict convexity of the master function.
Estimates for the momenta gaps of the quantum particles play a
role in the study of the thermodynamical limit,
see \cite{LL} and \cite{YY}. See e.g. \cite[Sect 4.3.2]{Ga} for
the exact analog of the estimates \eqref{momentgaps} for $\Sigma_0$ of type $A_n$.

It is believed \cite{IK} that quantum integrable systems governed by a strictly convex
master function always have a regularity constraint on the spectrum, although a conceptual
understanding is not known as far as we know.
We remark though that our derivation of the regularity constraint on the spectrum is in accordance to
this point of view. A conceptual understanding of the partly fermionic
nature of the quantum integrable system at hand is given in the
next section.


\section{Generalization of Girardeau's isomorphism}\label{sect3}

Let $C^\omega(V)$ be the space of complex valued, real analytic functions on
$V$, which we consider as a $W$-module with respect to the usual action \eqref{usual}.
Consider for $\lambda\in V_{\mathbb{C}}^*$ the space
\begin{equation}\label{Edef}
E(\lambda)=\{f\in C^\omega(V) \,\,\, | \,\,\,
p(\partial)f=p(\lambda)f\,\,\,\quad \forall\,p\in
S(V)_{\mathbb{C}}^{W_0} \},
\end{equation}
which is a $W$-submodule of $C^\omega(V)$.
We observed in Example \ref{k=nul} that
\begin{equation}\label{regularform}
E(\lambda)=\textup{BVP}_0(\lambda),\qquad \lambda\in
V_{\mathbb{C}}^*.
\end{equation}
In this section we give a convenient description of the
solution space $\hbox{BVP}_k(\lambda)^W$ of the boundary value
problem (Definition \ref{BVP}) in terms of the space of invariants in $E(\lambda)$
with respect to a $k$-dependent $W$-action by integral-reflection operators.
We will view this result as a natural generalization of Girardeau's \cite{Gi}
equivalence between the impenetrable quantum Bose-gas
and the free quantum Fermi-gas on the circle to arbitrary root
systems and to arbitrary multiplicity functions $k$.

We start by generalizing Girardeau's \cite{Gi} results on the impenetrable quantum
Bose-gas on the circle to arbitrary affine root systems.
Denote $E(\lambda)^{Q^\vee}$ for the subspace of $Q^\vee$-translation
invariant functions in $E(\lambda)$.

\begin{lem}\label{Qlem} For $\lambda\in V_{\mathbb{C}}^*$,
we have $E(\lambda)^{Q^\vee}=\{0\}$ unless $\lambda\in 2\pi i P$.
For $\lambda\in 2\pi i P$ the space $E(\lambda)^{Q^\vee}$ is spanned by
$e^\mu$ \textup{(}$\mu\in W_0\lambda$\textup{)}.
\end{lem}
\begin{proof}
By \cite{S}, a function $f\in E(\lambda)$ can be uniquely
expressed as $f(v)=\sum_{\mu\in W_0\lambda}p_\mu(v)e^{\mu(v)}$
where $p_\mu\in P(V)_{\mathbb{C}}$, see also Section \ref{sect7}.
Such a nonzero function $f$ is $Q^\vee$-translation
invariant iff
\begin{equation}\label{trp}
p_\mu(v+\gamma)e^{\mu(\gamma)}=p_\mu(v)
\end{equation}
for all $\mu\in
W_0\lambda$, $v\in V$ and $\gamma\in Q^\vee$.
This implies that $\lambda\in iV^*$ and that $p_\mu$ is bounded on $V$ for all $\mu\in W_0\lambda$. The
latter condition implies that $p_\mu$ is constant for all $\mu\in
W_0\lambda$. Returning to \eqref{trp} with $p_\mu\in\mathbb{C}$, the $Q^\vee$-translation
invariance of $f$ is equivalent to $\mu(Q^\vee)\subset 2\pi i\mathbb{Z}$
if $p_\mu\not=0$. Hence $E(\lambda)^{Q^\vee}=\{0\}$ unless
$\lambda\in 2\pi i P$, in which case $E(\lambda)^{Q^\vee}$ is
spanned by $e^\mu$ ($\mu\in W_0\lambda$).
\end{proof}

We denote $E(\lambda)^{-W}$ for the space of functions $f\in E(\lambda)$
satisfying $f(w^{-1}v)=(-1)^{l(w)}f(v)$ for all $w\in W$ and $v\in V$.
Since translations $\mu\in Q^\vee\subset W$ have even length,
$E(\lambda)^{-W}$ consists of $Q^\vee$-translation invariant functions.
In particular, $E(\lambda)^{-W}$ is
the solution space to the spectral problem for free fermionic
quantum integrable system on $V/Q^\vee$ associated to the Laplacean $\Delta$
on $V$.
\begin{cor}\label{fermecase}
Let $\lambda\in V_{\mathbb{C}}^*$. The space
$E(\lambda)^{-W}$ is zero-dimensional or one-dimensional.
It is one-dimensional iff $\lambda$ is a regular element from $2\pi i P$, in
which case $E(\lambda)^{-W}$ is spanned by
\begin{equation}\label{psiinftyextreme}
\psi_\lambda^\infty=\frac{1}{\#W_0}\prod_{\alpha\in\Sigma_0}\lambda(\alpha^\vee)^{-1}
\sum_{w\in
W_0}(-1)^{l(w)}e^{w\lambda}.
\end{equation}
\end{cor}
\begin{proof}
Let $\lambda\in 2\pi i P$ and $f=\sum_{\mu\in W_0\lambda}c_\mu e^{\mu}\in
E(\lambda)^{Q^\vee}$ with $c_\mu\in\mathbb{C}$, cf. Lemma \ref{Qlem}.
Then we have $f\in E(\lambda)^{-W}$ iff $c_{w\lambda}=(-1)^{l(w)}c_\lambda$ for all
$w\in W_0$. For singular $\lambda$ this implies $c_\mu=0$ for all $\mu\in W_0\lambda$.
For regular $\lambda$ we conclude that $f$ is a constant multiple of
$\psi_\lambda^\infty\in E(\lambda)^{-W}$.
\end{proof}
Following the analogy with Girardeau's \cite{Gi} analysis
of the impenetrable quantum Bose-gas on the circle,
we define now a linear map $G: C^\omega(V)\rightarrow C(V)^W$ by
\begin{equation}\label{G}
\bigl(Gf\bigr)(w^{-1}v):=f(v),\qquad w\in
W,\,\,v\in\overline{C_+}.
\end{equation}
The map $G$ is injective: for $g\in C(V)^W$ in the image of $G$,
the function $G^{-1}g$ is the unique analytic continuation of
$g|_{C_+}$ to $V$.

For $k\equiv \infty$ we interpret the boundary
conditions \eqref{jumpcond} as $f|_{V_{a}}\equiv 0$ for all $a\in \Sigma^+$.
The solution spaces $\textup{BVP}_\infty(\lambda)^W$ of the
associated boundary value problem
(see Definition \ref{BVP}) can now be analyzed as follows.

\begin{prop}\label{Girardeau}
For  $\lambda\in V_{\mathbb{C}}^*$ we have\\
{\bf (i)} The map $G$ restricts to a linear isomorphism
$G: E(\lambda)^{-W}\overset{\sim}{\longrightarrow}
\textup{BVP}_\infty(\lambda)^W$.\\
{\bf (ii)} The space $\textup{BVP}_\infty(\lambda)^W$ is
zero-dimensional or one-dimensional. It is one-dimensional iff
$\lambda$ is a regular element from $2\pi i P$. In that case
$\textup{BVP}_\infty(\lambda)^W$ is spanned by $\phi_\lambda^\infty:=
G(\psi_\lambda^\infty)$, which is the unique $W$-invariant
function satisfying
\[\phi_\lambda^\infty(v)=\frac{1}{\#W_0}\prod_{\alpha\in\Sigma_0^+}\lambda(\alpha^\vee)^{-1}
\sum_{w\in W_0}(-1)^{l(w)}e^{w\lambda(v)},\qquad v\in \overline{C_+}.
\]
\end{prop}
\begin{proof}
{\bf (i)} A function $f\in E(\lambda)^{-W}$ vanishes on the root
hyperplanes $V_a$ ($a\in\Sigma^+$), hence so does $g:=G(f)\in
C(V)^W$. The function $g$ furthermore satisfies the differential
equations \eqref{diffeqn}, hence $g\in \textup{BVP}_\infty(\lambda)^W$.

For $g\in \textup{BVP}_\infty(\lambda)^W$ we define $f=\widetilde{G}(g)\in
C(V)^{-W}$ by $f(w^{-1}v):=(-1)^{l(w)}g(v)$ for $w\in W$ and $v\in\overline{C_+}$.
This is well defined since $g$ vanishes on the root hyperplanes $V_a$ ($a\in\Sigma^+$).
Since $f$ is $W$-alternating we have
$f\in C^{1,(0)}(V)$. The function $f$ satisfies the differential
equations \eqref{diffeqn}, hence $f\in
\textup{BVP}_0(\lambda)^{-W}=E(\lambda)^{-W}$,
where the last equality follows from \eqref{regularform}.
The proof is now completed by observing that
$\widetilde{G}: \textup{BVP}_\infty(\lambda)^W\rightarrow
E(\lambda)^{-W}$ is the inverse of the map $G: E(\lambda)^{-W}\rightarrow
\textup{BVP}_\infty(\lambda)^W$.\\
{\bf (ii)} This follows from {\bf (i)} and Corollary
\ref{fermecase}.
\end{proof}
For root system $\Sigma_0$ of type $A$, Proposition
\ref{Girardeau} is due to Girardeau \cite{Gi}.

For the generalization of Proposition \ref{Girardeau} to arbitrary
multiplicity function $k$ it is convenient to reinterpret the space $E(\lambda)^{-W}$
as follows. Consider the integral operator
\begin{equation}\label{I}
\bigl(\mathcal{I}(a)f\bigr)(v)=\int_0^{a(v)}f(v-tDa^\vee)dt\qquad
(a\in\Sigma)
\end{equation}
as linear operator on $C(V)$. The integral operators
$\mathcal{I}(a)$ ($a\in I$) satisfy the braid relations of
$\Sigma$ as well as the quadratic relations $\mathcal{I}(a)^2=0$,
cf. e.g. \cite{Gu2}. In particular, given a reduced expression
$w=s_{i_1}s_{i_2}\cdots s_{i_{l(w)}}$ for $w\in W$, the operator
$Q_{\infty}(w):=\mathcal{I}(a_{i_1})\mathcal{I}(a_{i_2})\cdots \mathcal{I}(a_{i_{l(w)}})$
is well defined. Denote
\[E(\lambda)_{Q_\infty}^W:=\{f\in E(\lambda) \, | \,
Q_\infty(w)f=0,\quad \forall\,w\in W\setminus \{e\}\},
\]
where $e\in W$ is the unit element of $W$.
We now have the following simple observation.
\begin{lem}\label{simplelemma}
For $f\in C(V)$ and $b\in \Sigma$ we have $s_bf=-f$ if and only if
$\mathcal{I}(b)f=0$. In particular, $E(\lambda)^{-W}=E(\lambda)_{Q_\infty}^W$
for all $\lambda\in V_{\mathbb{C}}^*$.
\end{lem}
\begin{proof}
It is immediate that $\mathcal{I}(b)f=0$ if $s_bf=-f$. The converse follows from
the fact that
\begin{equation}
\label{Ianti}
\partial_{Db^\vee}\bigl(\mathcal{I}(b)f\bigr)=f+s_bf,
\end{equation}
cf. \cite[Lem. 2.1(iii)]{G}.
\end{proof}
By Lemma \ref{simplelemma}, Proposition
\ref{Girardeau}{\bf (i)} can be reformulated as the statement that
the map $G$ restricts to an isomorphism
\begin{equation}\label{Gisoref}
G: E(\lambda)_{Q_{\infty}}^W\overset{\sim}{\longrightarrow}
\textup{BVP}_\infty(\lambda)^W.
\end{equation}
The isomorphism \eqref{Gisoref} can now be generalized to
arbitrary multiplicity function $k$ as follows.
In the terminology of Gutkin \cite{Gu2},
the system of integral operators $\{k_b\mathcal{I}(b)\}_{b\in\Sigma^+}$ is an operator calculus with respect
to the affine Weyl group $W$ for arbitrary multiplicity function $k$.
This implies that the assignment $s_a\mapsto Q_{k,a}$ ($a\in I$),
with $Q_{k,a}$ the integral-reflection operators
\begin{equation}\label{Q}
\bigl(Q_{k,a}f\bigr)(v)=f(s_av)+k_a\bigl(\mathcal{I}(a)f\bigr)(v),\qquad
a\in \Sigma,\,\,f\in C(V),
\end{equation}
uniquely defines a $W$-action on $C(V)$,
cf. \cite{G}, \cite{Gu2} or Section \ref{sect5}. Accordingly, we
write
\[
Q_k(w):=Q_{k,a_{i_1}}Q_{k,a_{i_2}}\cdots Q_{k,a_{i_{r}}}
\]
for $w=s_{i_1}s_{i_2}\cdots s_{i_r}\in W$. Note that
\[Q_\infty(w)=\lim_{k\rightarrow\infty}k_w^{-1}Q_k(w),\qquad \forall\,w\in W,
\]
where $k_w:=k_{a_{i_1}}k_{a_{i_2}}\cdots k_{a_{i_r}}$ for
a reduced expression $w=s_{i_1}s_{i_2}\cdots s_{i_r}\in W$.
The generalization of \eqref{Gisoref} for arbitrary multiplicity function $k$
is now the statement that the map $G$ restricts to a linear isomorphism
\begin{equation}\label{Gkiso}
G: E(\lambda)_{Q_k}^W\overset{\sim}{\longrightarrow}
\textup{BVP}_k(\lambda)^W
\end{equation}
for arbitrary positive multiplicity function $k$, where $E(\lambda)_{Q_k}^W$
is the subspace of $Q_k(W)$-invariant functions in $E(\lambda)$. The proof of
\eqref{Gkiso} will be given in Section \ref{sect6}.

With the isomorphism \eqref{Gkiso}
at hand, Theorem \ref{main} is equivalent to the following theorem.
\begin{thm}\label{main2}
Let $\lambda\in V_{\mathbb{C}}^*$. The space $E(\lambda)_{Q_k}^{W}$
is one-dimensional or zero-dimensional. It is one-dimensional if and only if
$\lambda$ is a purely imaginary, regular solution of the Bethe ansatz equations
\eqref{BAE}. If $E(\lambda)_{Q_k}^{W}$ is one-dimensional then
\begin{equation}\label{BAF2}
\psi_\lambda^k(v)=\frac{1}{\#W_0}\sum_{w\in
W_0}\widetilde{c}_k(w\lambda)e^{w\lambda(v)},\qquad \forall v\in V
\end{equation}
is the unique function in $E(\lambda)_{Q_k}^{W}$ normalized by $\psi_\lambda^k(0)=1$.
\end{thm}
Theorem \ref{main2} is proved in Section \ref{sect8} under the assumption that $\lambda$ is regular.
The assertion that $\lambda$ is necessarily regular is proved in Section \ref{sect11}.

In order to reveal the full symmetry structures underlying the isomorphism \eqref{Gkiso}, we
will consider the upgrade of the map $G$ to a $k$-dependent linear
isomorphism $T_k$ of $C(V)$ which intertwines the
$Q_k(W)$-action with the usual $W$-action \eqref{usual}, and which acts as $G$ when
applied to $Q_k(W)$-invariant functions. The map which does the job is
Gutkin's \cite{G} propagation operator, defined by
$\bigl(T_kf\bigr)(w^{-1}v)=\bigl(Q_k(w)f\bigr)(v)$ for $w\in W$
and $v\in \overline{C_+}$ (see Section \ref{sect5} for details).
The propagation operator $T_k$ now restricts to an
isomorphism
\begin{equation}\label{Tisointro}
T_k: E(\lambda)\overset{\sim}{\longrightarrow}
\textup{BVP}_k(\lambda)
\end{equation}
for all $\lambda\in V_{\mathbb{C}}^*$ (cf. \cite{G} and Theorem
\ref{imageTlambda}), which implies \eqref{Gkiso} by restricting to
the subspaces of $W$-invariant functions.

We conclude this section by considering the limit
to the impenetrable case $k\equiv \infty$.
The Bethe ansatz equations \eqref{BAE} then reduce to
\[e^{w\lambda(\varphi^\vee)}=1\qquad \forall\,w\in W_0,
\]
which has $2\pi i P$ as purely imaginary solutions $\lambda$
(see \eqref{Palt}). Furthermore we have
\begin{equation}\label{muextremal}
\lim_{k\rightarrow\infty}\widehat{\mu}_k=2\pi\mu
\end{equation}
for $\mu\in P^{+}$, which follows by taking the limit
$k\rightarrow\infty$ in \eqref{momentgaps}.
For $\lambda=i\widehat{\mu}_k\in iV^*$ ($\mu\in P^{++}$) a regular
solution to the Bethe ansatz equation,
$\psi_\lambda^k\in E(\lambda)_{Q_k}^W$ (see \eqref{BAF2}) can alternatively be written as
\[\psi_\lambda^k=\frac{1}{\#W_0}\sum_{w\in
W_0}Q_k(w)(e^{\lambda}),
\]
see \cite{HO} or Section \ref{sect7}. It follows that
\begin{equation*}
\lim_{k\rightarrow\infty}k_{w_0}^{-1}\psi_{i\widehat{\mu}_k}^k=
\frac{1}{\#W_0}Q_{\infty}(w_0)(e^{2\pi i\mu})=\psi_{2\pi i\mu}^\infty
\end{equation*}
for $\mu\in P^{++}$, uniformly on compacta.
Pulling the limits through the map $G$, we obtain
\[\lim_{k\rightarrow \infty}k_{w_0}^{-1}\phi_{i\widehat{\mu}_k}^k=
\phi_{2\pi i\mu}^\infty\]
for $\mu\in P^{++}$, uniformly on compacta.


\section{Dunkl operators and Hecke algebras}\label{sect4}

It is well known that conserved integrals for quantum integrable systems of Calogero-Moser type
can be conveniently expressed in terms of Dunkl-type operators, which are explicit
commuting first-order differential-reflection operators, see e.g. \cite{He}, \cite{Du} and \cite{BFV}.
The Dunkl operators, together with the usual Weyl group action \eqref{usual}, form a faithful
representation of suitable degenerations of affine Hecke algebras, see \cite[Cor. 2.9]{Opd}. The exploration
of these structures has been instrumental in solving the corresponding quantum integrable systems.

In this section we derive the Dunkl-type operators and the underlying Hecke algebra structures
for the periodic quantum integrable systems with delta-potentials as introduced in
Section \ref{sect2}. We initially define the Dunkl operators as explicit differential-reflection operators on the space
$C^\infty(V_{reg})$ of smooth functions on $V_{reg}$. In Section \ref{sect6} we obtain the key result that
these Dunkl operators act on the solution space
$\textup{BVP}_k(\lambda)$ to the boundary value problem. Together with the
usual $W$-action \eqref{usual}, the space $\textup{BVP}_k(\lambda)$ then becomes
a module over the associated graded algebra
$H_k$ of Cherednik's \cite{C} (suitably filtered) degenerate double affine Hecke
algebra.

On the other hand, we will show in Section \ref{sect5} that the $W$-action $Q_k$ on $E(\lambda)$
together with the directional derivatives $\partial_v$ ($v\in V$)
makes $E(\lambda)$ into a $H_k$-module.
With these upgraded symmetry structures, Gutkin's propagation operator $T_k$ turns out to
yield an isomorphism
\[
T_k: E(\lambda)\overset{\sim}{\longrightarrow}
\textup{BVP}_k(\lambda)
\]
of $H_k$-modules for all $\lambda\in V_{\mathbb{C}}^*$. It is this particular isomorphism
which is explored in Section \ref{sect6}
to (re-)prove and clarify crucial results on the boundary value problem (see Definition \ref{BVP}),
as well as on the associated bosonic theory.

We denote $\chi: \mathbb{R}\setminus \{0\}\rightarrow \{0,1\}$
for the characteristic function of the interval $(-\infty,0)$, so
$\chi(x)=1$ if $x<0$ and $\chi(x)=0$ if $x>0$.
For $a\in\Sigma$ the function $\chi_a(v):=\chi(a(v))$ ($v\in V_{reg}$) defines a smooth function
on $V_{reg}$, which is constant on the alcoves $C$ of $V_{reg}$.
In fact, for $w\in W$ and $a\in \Sigma^+$ we have
\begin{equation}\label{chia}
\chi_a|_{w^{-1}C_+}\equiv
\begin{cases}
1\quad &\hbox{ if } wa\in \Sigma^-\\
0\quad &\hbox{ if } wa\in \Sigma^+,
\end{cases}
\end{equation}
hence $\chi_a$ is nonzero on a given alcove $w^{-1}C_+$ ($w\in W$) for only finitely
many positive roots $a\in \Sigma^+$. The Dunkl-type operators
\begin{equation}\label{Dunkl}
\mathcal{D}_v^k=\partial_v+\sum_{a\in\Sigma^+}k_aDa(v)\chi_a(\cdot)s_a\qquad (v\in V),
\end{equation}
thus define linear operators on $C^\infty(V_{reg})$, which depend linearly on $v\in V$.
For $f\in C^\infty(V_{reg})$ and $w\in W$  we have by \eqref{chia}
\begin{equation}\label{alcoveform}
\mathcal{D}_v^kf|_{w^{-1}C_+}=\Bigl(\partial_vf+\sum_{a\in\Sigma^+\cap
w^{-1}\Sigma^-}k_aDa(v)s_af\Bigr)\bigg|_{w^{-1}C_+}.
\end{equation}
In particular, for the fundamental alcove $C_+$ we simply have
\begin{equation}\label{alcovesimpleform}
\mathcal{D}_v^kf|_{C_+}=\partial_vf|_{C_+}.
\end{equation}
The Dunkl operators $\mathcal{D}_v^k$ \textup{(}$v\in V$\textup{)} and the $W$-action \eqref{usual}
on $C^\infty(V_{reg})$ satisfy the following fundamental commutation relations.
\begin{thm}\label{fundcomm}
{\bf (i)} We have the cross relation
\[s_a\mathcal{D}_v^k=\mathcal{D}_{s_{Da}v}^ks_a+k_aDa(v),\qquad v\in V,\,\, a\in I.\]

{\bf (ii)} The Dunkl operators $\mathcal{D}_v^k$
\textup{(}$v\in V$\textup{)} pair-wise commute.
\end{thm}
\begin{proof}
{\bf (i)} Fix $v\in V$ and $a\in I$.
By a direct computation we have
\[s_a\mathcal{D}_v^ks_a=\partial_{s_{Da}v}+\sum_{b\in
s_a\Sigma^+}k_bDb(s_{Da}v)\chi_b(\cdot)s_b.
\]
Since $s_a\Sigma^+=\bigl(\Sigma^+\setminus \{a\}\bigr)\cup \{-a\}$ we obtain
\begin{equation*}
\begin{split}
s_a\mathcal{D}_v^k&=\mathcal{D}_{s_{Da}v}^ks_a-k_aDa(s_{Da}v)\bigl(\chi_a(\cdot)+\chi_{-a}(\cdot)\bigr)\\
&=\mathcal{D}_{s_{Da}v}^ks_a+k_aDa(v),
\end{split}
\end{equation*}
which is the desired cross relation.

{\bf (ii)}  We derive the commutativity of the
Dunkl operators $\mathcal{D}_v^k$ ($v\in V$) as a direct consequence of
\eqref{alcovesimpleform} and the cross relation.
Let $f\in C^\infty(V_{reg})$ and $v,v^\prime\in V$. We show by
induction on the length $l(w)$ of $w\in W$ that
\begin{equation}\label{stepind}
\lbrack \mathcal{D}_v^k,\mathcal{D}_{v^\prime}^k\rbrack
f|_{w^{-1}C_+}=0.
\end{equation}
By \eqref{alcovesimpleform}, equation \eqref{stepind} is obviously valid for $w=e$ the unit element of $W$.
To prove the induction step, it suffices to show that
\begin{equation}\label{commutatorW}
s_a\lbrack \mathcal{D}_v^k,\mathcal{D}_{v^\prime}^k\rbrack=
\lbrack \mathcal{D}_{s_{Da}v}^k,\mathcal{D}_{s_{Da}v^\prime}^k\rbrack s_a
\end{equation}
for all $a\in I$. For the proof of \eqref{commutatorW}, first observe that
\begin{equation}\label{crossthrough}
s_a\mathcal{D}_v^k\mathcal{D}_{v^\prime}^k-\mathcal{D}_{s_{Da}v}^k\mathcal{D}_{s_{Da}v^\prime}^ks_a=
k_a\bigl(Da(v^\prime)\mathcal{D}_v^k+Da(v)\mathcal{D}_{v^\prime}^k-Da(v)Da(v^\prime)\mathcal{D}_{Da^\vee}^k\bigr)
\end{equation}
for all $a\in I$, which follows from applying the cross relation twice.
Now \eqref{commutatorW} follows from the fact that the right hand side of \eqref{crossthrough} is
symmetric in $v$ and $v^\prime$.
\end{proof}

By Theorem \ref{fundcomm}{\bf (ii)}, the assignment $v\mapsto \mathcal{D}_v^k$ uniquely extends to an algebra morphism
$S(V)_{\mathbb{C}}\rightarrow \hbox{End}(C^\infty(V_{reg}))$. We denote $p(\mathcal{D}^k)$ for the differential-reflection
operator on $C^\infty(V_{reg})$ associated to $p\in S(V)_{\mathbb{C}}$.

\begin{thm}\label{Hdefthm}
{\bf (i)}
There exists a unique complex unital associative algebra $H_k=H_k(\Sigma)$ satisfying
\begin{enumerate}
\item[{\bf (a)}] $H_k=S(V)_{\mathbb{C}}\otimes \mathbb{C}[W]$ as vector spaces, with $\mathbb{C}[W]$ the
group algebra of $W$.
\item[{\bf (b)}] The maps $p\mapsto p\otimes e$ and $w\mapsto 1\otimes w$, with $e\in W$ the unit element
of $W$, are algebra embeddings of $S(V)_{\mathbb{C}}$ and $\mathbb{C}[W]$ into $H_k$.
\item[{\bf (c)}] The cross relations
\[
s_a\cdot v-\bigl(s_{Da}v\bigr)\cdot s_a=k_aDa(v)
\]
holds in $H_k$ for $a\in I$ and $v\in V\subset S(V)_{\mathbb{C}}$. Here we have identified $S(V)_{\mathbb{C}}$
and $\mathbb{C}[W]$ with their images in $H_k$ through the algebra embeddings of {\bf (b)}.
\end{enumerate}
{\bf (ii)} The assignment $p\mapsto p(\mathcal{D}^k)$ \textup{(}$p\in S(V)_{\mathbb{C}}$\textup{)},
together with the $W$-action \eqref{usual},
defines a faithful representation $\pi_k: H_k\rightarrow
\textup{End}(C^\infty(V_{reg}))$.
\end{thm}
\begin{proof}
Suppose that $\sum_{w\in W}p_w(\mathcal{D}^k)w=0$ as an
endomorphism of $C^\infty(V_{reg})$,
with only finitely many $p_w\in S(V)_{\mathbb{C}}$'s non zero. We show that
all $p_w$'s are zero.
Equation \eqref{alcovesimpleform} implies
\begin{equation}\label{zerofund}
\sum_{w\in W}p_w(\partial)(wf)|_{C_+}\equiv 0,\qquad f\in C^\infty(V_{reg}).
\end{equation}
Applying \eqref{zerofund} to functions $f$ of the form $u^{-1}g$ with $u\in W$ and
with $g\in C^\infty(V_{reg})$ having support
in the fundamental alcove $C_+$, we conclude that $p_u(\partial)=0$ as a constant coefficient differential
operator on smooth functions in some open ball $D\subset C^+$, hence $p_u=0$.

The proof of the theorem is now standard: let $\widetilde{H}_k$ be the complex
unital associative algebra generated by $v\in V$
and $s_a$ ($a\in I$) with defining relations as in {\bf (b)} and {\bf (c)} (so the vectors $v\in V$ pair-wise commute, the $s_a$ ($a\in I$)
are involutions satisfying the Coxeter relations associated to $\Sigma$ and $I$,
and the generators satisfy the cross relations from {\bf (c)}).
By Theorem \ref{fundcomm} and by the paragraph preceding this theorem, the assignment $v\mapsto \mathcal{D}_v^k$,
together with the $W$-action \eqref{usual},
uniquely defines an algebra morphism $\pi_k: \widetilde{H}_k\rightarrow \hbox{End}(C^\infty(V_{reg}))$.
By the previous paragraph and by the cross relations in $\widetilde{H}_k$
it follows that $\pi_k$ is injective and that
$\widetilde{H}_k\simeq S(V)_{\mathbb{C}}\otimes \mathbb{C}[W]$ as vector spaces (the Poincar{\'e}-Birkhoff-Witt Theorem
for $\widetilde{H}_k$). Both statements of the theorem are now
immediately clear.
\end{proof}
We use the notation $M_{\pi_k}$ to indicate that a
subspace $M\subseteq C^\infty(V_{reg})$
is a $W$-submodule or $H_k$-submodule of $C^\infty(V_{reg})$
with respect to the $\pi_k$-action.

\begin{rema}
If the values $k_a$ of the multiplicity function $k$ are considered to be independent central variables in
the definition of $H_k$, then
$H_k$ is graded by imposing the degree of $w\in W$ to be zero and the degrees of $v\in V$
and $k_a$ to be one. As a graded algebra, $H_k$ is the associated graded algebra of Cherednik's \cite{C}
degenerate double affine Hecke algebra $\mathbb{H}_k$, considered as a filtered
algebra by the
same degree function (the only difference in the definition of $\mathbb{H}_k$
is the cross relation (see Theorem \ref{Hdefthm}{\bf (c)}), which now is of the form
\[
s_a\cdot v-\bigl(s_{a}v\bigr)\cdot s_a=k_aDa(v)
\]
for $a\in I$, where $S(V)$ is considered as a $W$-module algebra
with the action of $s_0$ defined by $s_{0}v=s_{\varphi}(v)+
2\|\varphi\|^{-2}\varphi(v)1\in S(V)$).
\end{rema}

\begin{lem}\label{center}
The center $Z(H_k)$ of $H_k$ contains $S(V)_{\mathbb{C}}^{W_0}$.
\end{lem}
\begin{proof}
Observe that the cross relations in $H_k$ (see Theorem \ref{Hdefthm}{\bf (c)}) imply
\begin{equation}\label{crossmore}
s_a\cdot p-\bigl(s_{Da}p\bigr)\cdot s_a=-k_a\frac{\bigl(s_{Da}p\bigr)-p}{Da^\vee}
\end{equation}
for $a\in I$ and $p\in S(V)_{\mathbb{C}}$. It follows from
\eqref{crossmore} that $S(V)_{\mathbb{C}}^{W_0}\subseteq Z(H_k)$.
\end{proof}
\begin{rema}\label{centerremark}
Observe that the subalgebra $H^{(0)}_k\subset H_k$ generated by $W_0$ and $S(V)_{\mathbb{C}}$
is isomorphic to the degenerate affine Hecke algebra (also known as the graded Hecke algebra),
see e.g. \cite{HO} and \cite{L}. By \cite[Prop. 4.5]{L} we have $Z(H^{(0)}_k)=S(V)_{\mathbb{C}}^{W_0}$.
\end{rema}

For trivial multiplicity parameters $k\equiv 0$, the operator $p(\mathcal{D}^0)$ ($p\in S(V)_{\mathbb{C}}$)
on $C^\infty(V_{reg})$ is the constant-coefficient differential operator $p(\partial)$ on $C^\infty(V_{reg})$.
We have the following striking fact when $p\in S(V)_{\mathbb{C}}$ is $W_0$-invariant.
\begin{cor}\label{regtriv}
For $p\in S(V)_{\mathbb{C}}^{W_0}$ we have $p(\mathcal{D}^k)=p(\partial)$ as operators on $C^\infty(V_{reg})$.
\end{cor}
\begin{proof}
Let $p\in S(V)_{\mathbb{C}}^{W_0}$ and $f\in C^\infty(V_{reg})$. By \eqref{alcovesimpleform} we have
$p(\mathcal{D}^k)f|_{C_+}=p(\partial)f|_{C_+}$. Let $w\in W$ and $v\in C_+$. By Lemma \ref{center}
applied twice (once with multiplicity function $k$, once with $k\equiv 0$), we have
\begin{equation*}
\begin{split}
\bigl(p(\mathcal{D}^k)f\bigr)(w^{-1}v)&=\bigl(p(\mathcal{D}^k)(wf)\bigr)(v)\\
&=\bigl(p(\partial)(wf)\bigr)(v)=\bigl(p(\partial)f\bigr)(w^{-1}v),
\end{split}
\end{equation*}
hence $p(\mathcal{D}^k)f=p(\partial)f$.
\end{proof}


\begin{rema}\label{finite}
The Dunkl operators $\mathcal{D}_v^k$, Theorem \ref{fundcomm}, Theorem \ref{Hdefthm} and Corollary \ref{regtriv}
have their obvious analogs in the context of finite root systems. In that case, the
Dunkl-type operators are
\[\partial_v+\sum_{\alpha\in\Sigma_0^+}k_\alpha\alpha(v)\chi_\alpha(\cdot)s_\alpha,\qquad
v\in V
\]
realizing, together with the $W_0$-action \eqref{usual}, an action of the degenerate affine Hecke algebra $H^{(0)}_k$
on the space of smooth functions on $V\setminus \bigcup_{\alpha\in\Sigma_0^+}V_\alpha$.
For classical root systems these operators were constructed using solutions of classical Yang-Baxter equations
and reflection equations in \cite{P}, \cite{MW} (type A) and \cite{KH}.
This construction fits into Cherednik's \cite{Ch0} general
framework relating root system analogs of $r$-matrices to (degenerate)
affine Hecke algebras and Dunkl operators.
\end{rema}


\section{Integral-reflection operators and the propagation operator}\label{sect5}

Heckman and Opdam \cite{HO} clarified the role of the degenerate affine Hecke
algebra $H^{(0)}_k$
in Gutkin's \cite{G} work when the underlying root system is finite.
It led to an explicit action of $H^{(0)}_k$ as directional derivatives
and integral-reflection operators. In this section we extend these results to
the present affine set-up. We show that Gutkin's \cite{G} propagation
operator intertwines this action with the action $\pi_k$ which is defined
in the previous section in
terms of Dunkl-type differential-reflection operators.

The integral-reflection operators $Q_{k,a}$ (see \eqref{Q}) for $a\in\Sigma$
are endomorphisms of $C(V)$ satisfying
\begin{equation}\label{QW}
wQ_{k,a}w^{-1}=Q_{k,w(a)},\qquad w\in W,\,\,a\in \Sigma
\end{equation}
with respect to the $W$-action \eqref{usual} on $C(V)$.
We furthermore have
\begin{equation}\label{trivial}
Q_{k,a}f|_{V_a}=f|_{V_a},\qquad a\in\Sigma.
\end{equation}
By \cite[Thm. 2.3]{G}, the assignment
\begin{equation}\label{Wpart}
s_a\mapsto Q_k(s_a):=Q_{k,a}\qquad (a\in I)
\end{equation}
extends to a representation $Q_k$ of $W$ on $C(V)$. In particular,
for $w\in W$ and any choice of decomposition
$w=s_{j_1}s_{j_2}\cdots s_{j_r}$ as a product of simple reflections
($j_l\in \{0,\ldots,n\}$), we have
\begin{equation}
Q_k(w)=Q_{k}(s_{j_1})Q_{k}(s_{j_2})\cdots Q_{k}(s_{j_r})
\end{equation}
as operators on $C(V)$.
\begin{defi}
Gutkin's \cite{G} propagation operator $T_k$ is the endomorphism of
$C(V)$ defined by
\begin{equation}\label{defT}
(T_kf)(w^{-1}v)=\bigl(Q_k(w)f)(v),\qquad v\in C_+,\quad w\in W.
\end{equation}
In particular, $T_0$ is the identity operator on $C(V)$.
\end{defi}
 A $W$-submodule $M\subseteq C(V)$ with respect to
the $Q_k$-action will be denoted by $M_{Q_k}$.
By construction the propagation operator $T_k: C(V)_Q\rightarrow C(V)_\pi$ is
$W$-equivariant. In fact, by \cite[Thm. 2.6]{G} $T_k$ is an isomorphism of
$W$-modules.

Observe that the operators $Q_k(w)$ ($w\in
W$) preserve the space $C^\infty(V)$ of complex valued, smooth functions
on $V$. The following result is the affine analog of
\cite[Thm. 2.1]{HO} and \cite[Cor. 2.3]{HO}.
\begin{thm}\label{Hdefthm2}
The assignment $v\mapsto \partial_v$ \textup{(}$v\in V$\textup{)}, together with the $W$-action
\eqref{Wpart} on $C^\infty(V)$, extends uniquely to a representation $Q_k: H_k\rightarrow \textup{End}(C^\infty(V))$.
\end{thm}
\begin{proof}
It suffices to verify the cross relations (see Theorem \ref{Hdefthm}{\bf (c)}), which
follow directly from \cite[Lem. 2.1]{G}.
\end{proof}
We will also use the notation $M_{Q_k}$ to indicate that a subspace $M\subseteq C^\infty(V)$ is a $H_k$-submodule
with respect to the $Q_k$-action. Observe that $C^\omega(V)_Q\subseteq C^\infty(V)_Q$
as $H_k$-submodule.

Consider the space $CB^\omega(V)$ of functions $f\in C(V)$ such that $f|_C$ is the restriction of
a \textup{(}necessarily unique\textup{)} analytic function on $V$ for all alcoves $C\in\mathcal{C}$
(cf. Remark \ref{ellarg}). Denote $C^{\omega,(k)}(V)$ for the
space of functions $f\in CB^\omega(V)$ satisfying
\begin{equation}\label{derivativejumpsk}
\partial^r_{Db^\vee}f\bigl(v+0Db^\vee\bigr)-\partial^r_{Db^\vee}f\bigl(v-0Db^\vee\bigr)=
\bigl(1-(-1)^r\bigr)k_b\partial^{r-1}_{Db^\vee}f\bigl(v+0Db^\vee\bigr)
\end{equation}
for $b\in\Sigma^+$, $v\in V_b$ sub-regular and $r\in\mathbb{Z}_{>0}$.
A function $f\in C^{\omega,(k)}(V)$ automatically
satisfies the jump conditions \eqref{derivativejumpsk}
for $b\in \Sigma^-$, $v\in V_b$ sub-regular and
$r\in\mathbb{Z}_{>0}$, hence the space $C^{\omega,(k)}(V)$ does not dependent on the choice
of positive roots $\Sigma^+$ in $\Sigma$. We thus can and will
interpret $C^{\omega,(k)}(V)_{\pi_k}$ and $CB^\omega(V)_{\pi_k}$ as $W$-submodules
of $C^\infty(V_{reg})_{\pi_k}$. Observe furthermore that $C^{\omega,(k)}(V)$
is a subspace of the space $C^{1,(k)}(V)$ used in the formulation of
the boundary value problems (see Proposition \ref{visual1} and Definition \ref{BVP}).

Observe that the propagation operator $T_k$ restricts to a linear map
\[
T_k: C^\omega(V)\rightarrow CB^\omega(V).
\]
We now obtain the following theorem.

\begin{thm}\label{imageT2}
{\bf (i)} $C^{\omega,(k)}(V)_{\pi_k}\subseteq
C^\infty(V_{reg})_{\pi_k}$ is a $H_k$-submodule.\\
{\bf (ii)} The propagation operator $T_k$ restricts to an isomorphism
\[
T_k: C^\omega(V)_{Q_k}\overset{\sim}{\longrightarrow}
C^{\omega,(k)}(V)_{\pi_k}
\]
of $H_k$-modules.
\end{thm}
\begin{proof}
We first show that $T_k$ restricts to a linear isomorphism
$T_k: C^\omega(V)\overset{\sim}{\longrightarrow} C^{\omega,(k)}(V)$.
For this we use the commutation relations
\begin{equation}\label{jumpH}
s_a\cdot\bigl(Da^\vee\bigr)^r-(-1)^r\bigl(Da^\vee\bigr)^r\cdot s_a=
\bigl(1-(-1)^r\bigr)k_a\bigl(Da^\vee\bigr)^{r-1},
\qquad a\in I,\,\, r\in \mathbb{Z}_{>0}
\end{equation}
in $H_k$, which follows from \eqref{crossmore}
applied to $p=(Da^\vee)^r\in S(V)_\mathbb{C}$.

Let $\phi\in C^\omega(V)$ and denote $f=T_k\phi\in CB^\omega(V)$.
We show that $f$ satisfies
the derivative jumps \eqref{derivativejumpsk} over
sub-regular $v\in V_b$ ($b\in\Sigma^+$) for all $r\in\mathbb{Z}_{>0}$.
In view of the $W$-equivariance of the propagation operator
$T_k$, it suffices to derive the derivative jumps for $f$ over sub-regular vectors
$v\in V_{a}\cap \overline{C_+}$ ($a\in I$).
Fix $a\in I$, $v\in V_a\cap \overline{C_+}$ sub-regular
and $r\in \mathbb{Z}_{>0}$.
For $\epsilon>0$ small we have $v+tDa^\vee=s_a(v-tDa^\vee)\in C_+$ for
$0<t<\epsilon$. Hence
\begin{equation}\label{aa1}
\partial^r_{Da^\vee}f(v+0Da^\vee)=\partial^r_{Da^\vee}\phi(v)=
Q_{k}(s_a)(\partial^r_{Da^\vee}\phi)(v),
\end{equation}
where the second equality follows from \eqref{trivial}. On the other hand,
\begin{equation}\label{aa2}
\partial^r_{Da^\vee}f(v-0Da^\vee)=(-1)^r\partial^r_{Da^\vee}(s_af)(v+0Da^\vee)
=(-1)^r\partial^r_{Da^\vee}(Q_{k}(s_a)\phi)(v).
\end{equation}
Combining \eqref{aa1} and \eqref{aa2} now yields
\begin{equation*}
\begin{split}
\partial^r_{Da^\vee}f(v+0Da^\vee)-\partial^r_{Da^\vee}f(v-0Da^\vee)&=
\bigl(\bigl(Q_k(s_a)\partial^r_{Da^\vee}-(-1)^r\partial^r_{Da^\vee}Q_k(s_a)\bigr)\phi\bigr)(v)\\
&=\bigl(1-(-1)^r\bigr)k_a\partial^{r-1}_{Da^\vee}\phi(v)\\
&=\bigl(1-(-1)^r\bigr)k_a\partial^{r-1}_{Da^\vee}f(v+0Da^\vee),
\end{split}
\end{equation*}
where the second equality follows from (the $Q_k$-image of)
\eqref{jumpH}. Thus $f\in C^{\omega,(k)}(V)$.

The map $T_k: C^\omega(V)\rightarrow C^{\omega,(k)}(V)$
is clearly injective. We now proceed to prove surjectivity.
Let $f\in C^{\omega,(k)}(V)$ and denote $\psi$
for the unique analytic function on $V$ satisfying $\psi|_{C_+}=f|_{C_+}$. The function
$g:=f-T_k\psi\in C^{\omega,(k)}(V)$ satisfies
$g|_{\overline{C_+}}\equiv 0$. Combined with the continuity of $g$ and the
derivative jump conditions \eqref{derivativejumpsk} for $g$, we obtain
\[\bigl(\partial_{Da^\vee}^rg\bigr)(v-0Da^\vee)=0
\]
for $r\in\mathbb{Z}_{\geq 0}$, $a\in I$ and $v\in V_a\cap \overline{C_+}$
sub-regular. Since $g|_C$ ($C\in\mathcal{C}$) has an extension to an analytic function on the
whole Euclidean space $V$, we conclude that $g|_{\overline{C}}\equiv 0$ for the
neighboring alcoves $C=s_aC_+$ ($a\in I$) of $C_+$. Continuing
inductively we conclude that $g\equiv 0$ on $V$, hence
$f=T_k\psi$.

It remains to show that the isomorphism
\[T_k: C^\omega(V)_{Q_k}\overset{\sim}{\longrightarrow} C^{\omega,(k)}(V)_{\pi_k}
\]
of $W$-modules is in fact an isomorpism of $H_k$-modules. For this it suffices to show that
\begin{equation}\label{ip}
T_k\bigl(\partial_vf\bigr)|_{V_{reg}}=\mathcal{D}_v^k(T_kf|_{V_{reg}})
\end{equation}
for $v\in V$ and $f\in C^\omega(V)$.
To prove \eqref{ip} we use the commutation relation
\begin{equation}\label{crossextended}
w\cdot v=\bigl((Dw)v\bigr)\cdot w+\sum_{a\in\Sigma^+\cap w^{-1}\Sigma^-}k_aDa(v)ws_a
\end{equation}
in $H_k$, which can be easily proved by induction on the length $l(w)$ of $w\in W$ using the
cross relations in $H_k$ (see Theorem \ref{Hdefthm}{\bf (c)}). Fix $w\in W$ and $v^\prime\in C_+$. By \eqref{crossextended}
and Theorem \ref{Hdefthm2} we have
\begin{equation*}
\begin{split}
T_k(\partial_vf)(w^{-1}v^\prime)&=Q_k(w)(\partial_vf)(v^\prime)\\
&=\partial_{(Dw)v}(Q_k(w)f)(v^\prime)+
\sum_{a\in \Sigma^+\cap w^{-1}\Sigma^-}k_aDa(v)Q_k(ws_a)f(v^\prime)\\
&=\partial_v(T_kf)(w^{-1}v^\prime)+\sum_{a\in\Sigma^+\cap
w^{-1}\Sigma^-}k_aDa(v)T_kf(s_aw^{-1}v^\prime)\\
&=\mathcal{D}_v^k(T_kf)(w^{-1}v^\prime),
\end{split}
\end{equation*}
where the last equality follows from \eqref{alcoveform}.
\end{proof}

\begin{rema}
The assertion \cite[Thm. 2.7]{G} that, in Gutkin's notation, the propagation operator $T_k$ is 
an automorphism of the $W$-module $CB^\infty$, seems to be incorrect. In fact, the integral-operators
$\mathcal{I}(a)$ ($a\in \Sigma$) do not preserve $CB^\infty$, contrary
to the claim in the proof of \cite[Thm. 2.7]{G}. In \cite{G}, this result is used to
link $\textup{BVP}_k(\lambda)$ to $E(\lambda)$ (see \eqref{Edef}). We will show in Section \ref{sect6} that 
Theorem \ref{imageT2}{\bf (ii)} suffices to provide this link. 
\end{rema}


\begin{rema}
Theorem \ref{imageT2} has an obvious analog in the context
of finite root systems (compare with Remark \ref{finite}).
In the case of finite root system of type A, the intertwining properties
of the propagation operator with respect to the degenerate affine Hecke algebra
actions were considered in \cite{H} and the normal derivative jump
conditions of higher order were considered in \cite{Gu}.

\end{rema}
\begin{cor}
Fix $v\in V$. The Dunkl operator $\mathcal{D}_v^k$ is a linear
operator on $C^{\omega,(k)}(V)$ satisfying
$\mathcal{D}_v^k\bigl(T_kf\bigr)=T_k\bigl(\partial_vf\bigr)$
for all $f\in C^\omega(V)$.
\end{cor}

In the following proposition we relate the
Dunkl operators $\mathcal{D}^k_v$ to the
quantum Hamiltonian $\mathcal{H}_k$ (see \eqref{Hamiltonian} and \eqref{distribution}).
Recall that $p_2(\partial)=\Delta$ for the $W_0$-invariant polynomial
$p_2=\|\cdot\|^2$ on $V^*$.
\begin{prop}\label{visual2}
For $f\in C^{\omega,(k)}(V)$ we have
\begin{equation}\label{Hformu}
-p_2(\mathcal{D}^k)f=\mathcal{H}_kf
\end{equation}
as distributions on $V$.
\end{prop}
\begin{proof}
Fix $f\in C^{\omega,(k)}(V)$, then $p_2(\mathcal{D}^k)f\in C^{\omega,(k)}(V)\subseteq C(V)$
and $p_2(\mathcal{D}^k)f|_{V_{reg}}=\Delta f|_{V_{reg}}$ by Corollary \ref{regtriv}. Furthermore, $f$
satisfies the first order normal derivative jumps \eqref{jumpcond}
over the affine hyperplanes $V_a$ ($a\in\Sigma^+$).
The identity \eqref{Hformu} then follows
from a standard argument using Green's identity,
cf. (the proof of) Proposition \ref{visual1}.
\end{proof}

By Proposition \ref{visual2} it is justified to interpret
the quantum Hamiltonian $\mathcal{H}_k$ on $C^{\omega,(k)}(V)$ as the
operator $-p_2(\mathcal{D}^k)$ on $C^{\omega,(k)}(V)$.
The complete integrability of the quantum system
is then directly reflected by the commutativity of the Dunkl
operators $\mathcal{D}_v^k$ ($v\in V$). More precisely, the space $C^{\omega,(k)}(V)_\pi^W$
serves as an algebraic model for the Hilbert space of quantum states
associated to the bosonic quantum system on $V/Q^\vee$ with Hamiltonian
$\mathcal{H}_k=-p_2(\mathcal{D}^k)$. The pair-wise commuting operators $p(\mathcal{D}^k)$ ($p\in S(V)_{\mathbb{C}}^{W_0}$)
on $C^{\omega,(k)}(V)_\pi^W$ are the corresponding quantum conserved integrals.


\section{The boundary value problem revisited}\label{sect6}

The operators $p(\mathcal{D}^k)$ ($p\in S(V)_{\mathbb{C}}^{W_0}$)
on $C^{\omega,(k)}(V)$ satisfy
\[p(\mathcal{D}^k)f|_{V_{reg}}=
p(\partial)f|_{V_{reg}},\qquad f\in C^{\omega,(k)}(V)
\]
by Corollary \ref{regtriv}. This key observation leads to an explicit connection
between the spectral problem of the operators $p(\mathcal{D}^k)$ ($p\in S(V)_{\mathbb{C}}^{W_0}$)
and the boundary value problem as formulated in Definition \ref{BVP}. We will first do the
analysis for the spectral problem of the quantum Hamiltonian $\mathcal{H}_k$
(defined by \eqref{Hamiltonian} and \eqref{distribution}).

For $E\in\mathbb{C}$ we write $\mathcal{E}(E)$ for the space of
functions $f\in C^{\omega}(V)$ satisfying $\Delta f=-Ef$ on $V$ (cf. Example \ref{k=nul}).
By Lemma \ref{center}, $\mathcal{E}(E)_{Q_k}\subseteq C^\omega(V)_{Q_k}$ is
a $H_k$-submodule.
 Denote $\mathcal{E}_k(E)$ for the space of functions
$f\in CB^\omega(V)$ satisfying $\mathcal{H}_kf=Ef$ as distributions on $V$.

\begin{thm}\label{regularize}
Fix $E\in\mathbb{C}$.\\
{\bf (i)}  The space $\mathcal{E}_k(E)$ is the $H_k$-submodule of $C^{\omega,(k)}(V)_{\pi_k}$
consisting of eigenfunctions of $p_2(\mathcal{D}^k)$
with eigenvalue $-E$.\\
{\bf (ii)} The propagation operator $T_k$ restricts to an isomorphism
\[T_k: \mathcal{E}(E)_{Q_k}\overset{\sim}{\longrightarrow} \mathcal{E}_k(E)_{\pi_k}
\]
of $H_k$-modules.
\end{thm}
\begin{proof}
{\bf (i)} We first show that $\mathcal{E}_k(E)\subseteq C^{\omega,(k)}(V)$.
Fix $f\in \mathcal{E}_k(E)$. By Proposition \ref{visual1},
$f\in C^{1,(k)}(V)\cap CB^\omega(V)$ and $\Delta f|_{V_{reg}}=-Ef|_{V_{reg}}$.
Let $\psi$ be the unique analytic function on $V$ satisfying
$\psi|_{C_+}=f|_{C_+}$, then $\psi\in\mathcal{E}(E)$. By Theorem \ref{imageT2} and Corollary
\ref{regtriv} we conclude that
$T_k\psi\in C^{\omega,(k)}(V)$ and $\Delta(T_k\psi)|_{V_{reg}}=-E(T_k\psi)|_{V_{reg}}$. Hence
\[
g:=f-T_k\psi\in C^{1,(k)}(V)\cap CB^\omega(V)
\]
satisfies $\Delta g|_{V_{reg}}=-Eg|_{V_{reg}}$ and has the additional property that
$g|_{\overline{C_+}}\equiv 0$. Fix $v\in V_a\cap \overline{C_+}$ ($a\in I$) sub-regular.
The nontrivial normal derivative jump condition \eqref{jumpcond} for $g$
at $v$ trivializes since $g|_{\overline{C_+}}\equiv 0$, hence $g$ is continuously differentiable
in an open neighborhood $U$ of $v$.
It follows that $g|_U$ is a distribution solution of the
(hypo)elliptic constant coefficient
differential operator $\Delta+E$ (cf. Example \ref{k=nul}),
hence $g|_U$ is smooth. Since $g|_{\overline{C_+}}\equiv 0$,
we conclude that
\[\partial_{Da^\vee}^rg(v-0Da^\vee)=\partial_{Da^\vee}^rg(v+0Da^\vee)=0,\qquad r\in\mathbb{Z}_{\geq 0}.
\]
As in the proof of Theorem \ref{imageT2} we conclude that
$g|_{\overline{s_aC_+}}\equiv 0$ for $a\in I$. Continuing inductively,
we conclude that $g\equiv 0$ on $V$. Hence $f=T_k\psi\in
C^{\omega,(k)}(V)$.

By Proposition \ref{visual2} we conclude that
\begin{equation}\label{centerform}
\mathcal{E}_k(E)=\{f\in C^{\omega,(k)}(V) \, | \, p_2(\mathcal{D}^k)f=-Ef
\}.
\end{equation}
Furthermore $p_2(\mathcal{D}^k)=\pi_k(p_2)$, hence Lemma
\ref{center} implies that
$\mathcal{E}_k(E)_{\pi_k}\subseteq C^{\omega,(k)}(V)_{\pi_k}$ is a $H_k$-submodule.\\
{\bf (ii)} This follows directly from Theorem \ref{imageT2},
\eqref{centerform} and the fact that $Q_k(p_2)=p_2(\partial)=\Delta$.

\end{proof}

We now extend these results to the solution spaces $\textup{BVP}_k(\lambda)$
of the boundary value problem (Definition \ref{BVP}).
For a $H_k$-module $M$ and $\lambda\in V_{\mathbb{C}}^*$ we define
\begin{equation}
M_\lambda:=\{m\in M \,\,\, | \,\,\, p\cdot m=
p(\lambda)m\quad \forall\, p\in S(V)_{\mathbb{C}}^{W_0}
\},
\end{equation}
which is a $H_k$-submodule of $M$ in view of Lemma \ref{center}.
By Remark \ref{centerremark} the module
$M_\lambda$ consists of the vectors $m\in M$ transforming
according to the central character
$\lambda\in V_{\mathbb{C}}^*$ for the action of
the center of the degenerate
affine Hecke algebra $H_k^{(0)}\subseteq H_k$.

\begin{cor}\label{regularizecor}
Let $\lambda\in V_{\mathbb{C}}^*$.
The space $\textup{BVP}_k(\lambda)$ is the $H_k$-submodule $C^{\omega,(k)}(V)_{\pi_k,\lambda}$
of $C^{\omega,(k)}(V)_{\pi_k}$.
\end{cor}
\begin{proof}
By Corollary \ref{regtriv} and Theorem \ref{imageT2} we have
\begin{equation}\label{secondk}
C^{\omega,(k)}(V)_{\pi_k,\lambda}=\{ f\in C^{\omega,(k)}(V) \,\,\, | \,\,\,
p(\partial)f|_{V_{reg}}=p(\lambda)f|_{V_{reg}} \quad \forall \, p\in
S(V)_{\mathbb{C}}^{W_0} \},
\end{equation}
hence $C^{\omega,(k)}(V)_{\pi_k,\lambda}\subseteq
\textup{BVP}_k(\lambda)$. By Proposition \ref{visual1} and Remark
\ref{ellarg} we have
\[\textup{BVP}_k(\lambda)\subseteq
\mathcal{E}_k(-p_2(\lambda)).
\]
Theorem \ref{regularize} and \eqref{secondk} now imply that
$\textup{BVP}_k(\lambda)\subseteq C^{\omega,(k)}(V)_{\pi_k,\lambda}$.
\end{proof}

\begin{thm}\label{imageTlambda}
Let $\lambda\in V_{\mathbb{C}}^*$.\\
{\bf (i)} The propagation operator $T_k$ restricts to an isomorphism
$T_k: E(\lambda)_{Q_k}\overset{\sim}{\longrightarrow}
\textup{BVP}_k(\lambda)_{\pi_k}$ of left $H_k$-modules.\\
{\bf (ii)} The map $G$ \eqref{G} restricts to an isomorphism
$G: E(\lambda)_{Q_k}^{W}\overset{\sim}{\longrightarrow}
\textup{BVP}_k(\lambda)_{\pi_k}^W$.
\end{thm}
\begin{proof}
{\bf (i)} The restriction of the propagation operator $T_k$ to
the $H_k$-module $E(\lambda)_{Q_k}=C^\omega(V)_{Q_k,\lambda}$ defines
an isomorphism
\[T_k: E(\lambda)_{Q_k}\overset{\sim}{\longrightarrow}
C^{\omega,(k)}(V)_{\pi_k,\lambda}
\]
of $H_k$-modules in view of Theorem \ref{imageT2}.
Corollary \ref{regularizecor} now completes
the proof.

{\bf (ii)} This follows from {\bf (i)} and from the fact that
the propagation map $T_k$ acts on $Q_k(W)$-invariant functions
in the same way as the map $G$ \eqref{G}.
\end{proof}

As observed in Section \ref{sect3}, Theorem \ref{imageTlambda}{\bf (ii)}
can be used to reformulate the main results on the solution
space $\hbox{BVP}_k(\lambda)_\pi^W$ (see Theorem \ref{main}) to the boundary value problem
in terms of the space of invariants $E(\lambda)_Q^W$, where
$E(\lambda)$ now is the solution space to the boundary value problem with
zero normal derivative jumps over sub-regular vectors. Theorem \ref{main2}
is the resulting reformulation of Theorem \ref{main}. In order to prove
Theorem \ref{main2} we analyze the space $E(\lambda)_Q^W$ in detail
in the following sections.


\section{Invariants in $E(\lambda)$}\label{sect7}

In this section we analyze the sub-space $E(\lambda)_Q^{W_0}$ of $W_0$-invariants of $E(\lambda)_Q$.
First we recall some well known properties of the space $E(\lambda)$ from \cite{S} and \cite{HO}.
For technical purposes it is convenient to introduce the following terminology.
\begin{defi}\label{Jstandard}
Let $J$ be a subset of the simple roots $I_0$. The spectral parameter
$\lambda\in V_{\mathbb{C}}^*$ is called $J$-standard
if $\lambda\in V^*\oplus i\overline{V_+^*}$ and if
the isotropic sub-group of $\lambda$ in $W_0$ is the standard parabolic
sub-group $W_{0,J}$ generated by the simple reflections $s_\alpha$
\textup{(}$\alpha\in J$\textup{)}.
\end{defi}

\begin{lem}\label{Jstandardlem}
Let $\lambda\in V_{\mathbb{C}}^*$. The $W_0$-orbit of $\lambda$ contains a $J$-standard spectral parameter
for some subset $J\subseteq I_0$.
\end{lem}
\begin{proof}
Taking a $W_0$-translate of $\lambda$ we may assume that $\lambda=\mu+i\nu$
with $\mu\in V^*$ and $\nu\in \overline{V_+^*}$. The isotropy group of $\nu$
in $W_0$ is a standard parabolic sub-group $W_{0,K}\subset W_0$ for some subset
$K\subseteq I_0$. Write $V^*=V_{K}^*\oplus (V_K^*)^\perp$
with $V_K^*=\hbox{span}_{\mathbb{R}}\{\alpha\, | \, \alpha\in K \}$ and $(V_K^*)^\perp$
its orthocomplement in $V^*$. Set
\[V_{K,+}^*=\{\xi\in V_K^* \, | \, \xi(\alpha^\vee)>0 \quad \forall\,\alpha\in K\},
\]
which we view as the fundamental chamber for the action of the
standard parabolic sub-group $W_{0,K}$
on $V_K^*$. Taking a $W_{0,K}$-translate of $\lambda$ we may
assume that $\lambda=\mu+\mu^\prime+i\nu$ with $\mu\in \overline{V_{K,+}^*}$, $\mu^\prime\in (V_{K}^*)^\perp$,
and $\nu\in \overline{V_+^*}$ as before. The isotropy sub-group of $\lambda$
in $W_0$ then equals the isotropy sub-group of $\mu$ in $W_{0,K}$, which is a standard
parabolic sub-group $W_{0,J}$ for some subset $J\subseteq K$ since $\mu\in \overline{V_{K,+}^*}$.
\end{proof}

Observe that a $J$-standard spectral parameter $\lambda$ is regular if and only if $J=\emptyset$.
Note furthermore that the module $E(\lambda)$ ($\lambda\in V_{\mathbb{C}}^*$) only depends on the orbit $W_0\lambda$.
When analyzing the module $E(\lambda)$, we thus may assume without loss of generality that $\lambda$
is $J$-standard for some subset $J\subseteq I_0$. In particular, we will now assume this condition
for the remainder of this section.

For $j\in\mathbb{Z}_{\geq 0}$ we
denote $P^{(j)}(V)_{\mathbb{C}}$ (respectively $P^{(\leq
j)}(V)_{\mathbb{C}}$) for the homogeneous polynomials $p\in P(V)_{\mathbb{C}}$
of degree $j$ (respectively the polynomials $p\in P(V)_{\mathbb{C}}$
of degree $\leq j$). The $W_0$-action \eqref{usual} on $P(V)_{\mathbb{C}}$
respects the natural grading
$P(V)_{\mathbb{C}}=\bigoplus_{j=0}^{\infty} P^{(j)}(V)_{\mathbb{C}}$.
Furthermore,
\[E_J(0)=\{f\in P(V)_{\mathbb{C}} \, | \,
p(\partial)f=p(0)f\quad\forall\, p\in S(V)^{W_{0,J}}\}
\]
is a graded $W_{0,J}$-submodule of $P(V)_{\mathbb{C}}$, isomorphic to the
regular representation of $W_{0,J}$ (see e.g. \cite[Thm. 1.2]{S} and
references therein). We write
$E_J^{(j)}(0)=E_J(0)\cap P^{(j)}(V)_{\mathbb{C}}$ and
$E_J^{(\leq j)}(0)=E_J(0)\cap P^{(\leq j)}(V)_{\mathbb{C}}$.

Denote by $W_0^J$ the minimal coset representatives of $W_0/W_{0,J}$.
Steinberg \cite{S} established the decomposition
\begin{equation}
\label{Steinberg}
E(\lambda)=\bigoplus_{u\in W_0^J}u\bigl(E_J(0)e^\lambda\bigr).
\end{equation}
It follows from \eqref{Steinberg} that $E(\lambda)$, viewed as a $W_0$-module by
the action \eqref{usual}, is isomorphic to the regular representation of $W_0$.
Furthermore, we have
$E(\lambda)=\bigoplus_{j=0}^\infty E^{(j)}(\lambda)$ with $E^{(j)}(\lambda)$
the $W_0$-submodule
\[E^{(j)}(\lambda)=\bigoplus_{u\in
W_0^J}u\bigl(E_J^{(j)}(0)e^\lambda\bigr).
\]
We denote $E^{(\leq j)}(\lambda)=\bigoplus_{r=0}^jE^{(r)}(\lambda)$.

Representations of the finite group $W_0$ do not admit nontrivial
continuous deformations, hence $E(\lambda)_{Q}$
is isomorphic to the regular representation of $W_0$ for arbitrary
multiplicity function $k$.
In particular, $E(\lambda)_Q^{W_0}$ is one-dimensional
for all spectral values $\lambda\in V_{\mathbb{C}}^*$. In fact, by
\eqref{trivial} the function
\begin{equation}\label{Qpsi}
\psi_\lambda^k=\frac{1}{\#W_0}\sum_{w\in W_0}Q_k(w)e^\lambda
\end{equation}
satisfies $\psi_\lambda^k(0)=1$ and spans $E(\lambda)_Q^{W_0}$.
On the other hand, by \eqref{Steinberg} there exist unique polynomials
$p^\lambda_u\in E_J(0)$ ($u\in W_0^J$) such that
\begin{equation}\label{exp}
\psi_\lambda^k(v)=\sum_{u\in
W_0^J}p^\lambda_u(u^{-1}v)e^{u\lambda(v)},
\qquad v\in V.
\end{equation}

By \eqref{Steinberg} we have
\begin{equation}\label{Steinbergregular}
E(\lambda)=\bigoplus_{w\in W_0}\mathbb{C}e^{w\lambda},
\qquad \lambda\in V_{\mathbb{C}}^*\,\,\hbox{regular},
\end{equation}
so the polynomials $p^\lambda_w$ ($w\in W_0$) are constants for regular $\lambda$.
In fact, from e.g. \cite{Ga} and \cite[Sect. 2]{HO} we have
\begin{equation}\label{regular}
\psi_\lambda^k(v)=\frac{1}{\#W_0}\sum_{w\in
W_0}\widetilde{c}_k(w\lambda)e^{w\lambda(v)},\qquad \lambda\in
V_{\mathbb{C}}^* \,\,\, \hbox{regular},
\end{equation}
where the $c$-function $\widetilde{c}_k$ is given by \eqref{cfunction}.
In the remainder of the paper it will actually be more
convenient to work with the regularized
$c$-function
\begin{equation}\label{regularc}
c_k(\mu):=\prod_{\stackrel{\alpha\in\Sigma_0^+}{\mu(\alpha^\vee)\not=0}}
\frac{\mu(\alpha^\vee)+k_\alpha}{\mu(\alpha^\vee)},\qquad \mu\in
V_{\mathbb{C}}^*
\end{equation}
which is equal to $\widetilde{c}_k(\mu)$ for regular $\mu$.
We can then write
\[p^\lambda_w=\frac{1}{\#W_0}c_k(w\lambda),\qquad
\lambda\in V_{\mathbb{C}}^*\,\,\,\hbox{regular}.
\]
For singular $\lambda$ an explicit expression for $p^\lambda_u\in E_J(0)$
($u\in W_0^J$) is not known. For our purposes it suffices to have
explicit expressions for the highest and the next to highest
homogeneous components of $p^\lambda_u$, which we will
now proceed to derive.

We denote $\Sigma_0^J\subseteq \Sigma_0$ for the parabolic root
sub-system associated to the simple roots $J\subseteq I_0$.
We write $N_J$ for the cardinality of the corresponding set $\Sigma_0^{J,+}:=\Sigma_0^J\cap \Sigma_0^+$ of
positive roots in $\Sigma_0^J$ and
\[\delta_J=\frac{1}{2}\sum_{\alpha\in\Sigma_0^{J,+}}\alpha\in V^*.
\]
Recall that the minimal coset representatives $W_0^J$ of $W_0/W_0^J$ can be characterized by
\[W_0^J=\{u\in W_0 \, | \, u(\Sigma_0^{J,+})\subseteq \Sigma_0^+ \}.
\]
The following lemma now gives a derivational expression for $p^\lambda_u$ ($u\in W_0^J$).
\begin{lem}\label{differformula}
Let $\lambda\in V_{\mathbb{C}}^*$ be $J$-standard. For $u\in W_0^J$ we have
\[
p^\lambda_u=K_J^{-1}
\frac{d^{N_J}}{dt^{N_J}}\bigg|_{t=0}\left(\sum_{v\in
W_{0,J}}d_u(t)e_{uv}(t)(-1)^{l(v)}e^{tv\delta_J}\right)
\]
with coefficients
\[
d_u(t)=\prod_{\alpha\in\Sigma_0^+\setminus u(\Sigma_0^{J,+})}\bigl(u\delta_J(\alpha^\vee)t+u\lambda(\alpha^\vee)\bigr)^{-1},\qquad
e_{uv}(t)=\prod_{\alpha\in\Sigma_0^+}(uv\delta_J(\alpha^\vee)t+u\lambda(\alpha^\vee)+k_\alpha)
\]
and with strictly positive constant $K_J=N_J!\#W_0\prod_{\alpha\in\Sigma_0^{J,+}}\delta_J(\alpha^\vee)$.
\end{lem}
\begin{proof}
By \eqref{Qpsi}, $\psi_\mu^k(v^\prime)$ ($v^\prime\in V$) depends analytically on the spectral parameter $\mu\in V_{\mathbb{C}}^*$.
In particular, $\psi_{\lambda_t}^k(v^\prime)$ with $\lambda_t:=\lambda+t\delta_J\in V_{\mathbb{C}}^*$
depends analytically on $t\in\mathbb{C}$, and we have the (point-wise) limit
\begin{equation}\label{limpsi}
\lim_{t\rightarrow 0}\psi_{\lambda_t}^k=\psi_\lambda^k.
\end{equation}
For $\epsilon>0$ we write
\[
U_\epsilon^0=\{t\in\mathbb{C} \, | \, 0<|t|<\epsilon\},\qquad
U_\epsilon=\{t\in \mathbb{C}\, | \, |t|<\epsilon \}.
\]
There exists an $\epsilon>0$ such that $\lambda_t$ is regular for $t\in U_\epsilon^0$,
hence
\[
\psi_{\lambda_t}^k=\frac{1}{\#W_0}\sum_{w\in
W_0}\left(\prod_{\alpha\in\Sigma_0^+}\frac{w\lambda_t(\alpha^\vee)+k_\alpha}{w\lambda_t(\alpha^\vee)}
\right)e^{w\lambda_t},\qquad t\in U_\epsilon^0
\]
by \eqref{regular}. Splitting the sum into a double sum $w=uv$ with $u\in W_0^J$
and $v\in W_{0,J}$ and using
\begin{equation*}
\begin{split}
\prod_{\alpha\in\Sigma_0^+}uv\lambda_t(\alpha^\vee)&=(-1)^{l(u)+l(v)}t^{N_J}\prod_{\alpha\in\Sigma_0^{J,+}}\delta_J(\alpha^\vee)
\prod_{\beta\in\Sigma_0^+\setminus\Sigma_0^{J,+}}\lambda_t(\beta^\vee)\\
&=(-1)^{l(v)}t^{N_J}\prod_{\alpha\in\Sigma_{0}^{J,+}}\delta_J(\alpha^\vee)
\prod_{\beta\in\Sigma_0^+\setminus
u(\Sigma_0^{J,+})}u\lambda_t(\beta^\vee),
\end{split}
\end{equation*}
we obtain
\begin{equation}\label{limpsi2}
t^{N_J}\psi_{\lambda_t}^k=K_J^{-1}N_J!
\sum_{u\in W_0^J}\sum_{v\in W_{0,J}}d_u(t)e_{uv}(t)(-1)^{l(v)}e^{tuv\delta_J+u\lambda}
\end{equation}
as analytic functions in $t\in U_\epsilon$ (note that $d_u(t)$ is analytic at
$t\in U_\epsilon$). By \eqref{limpsi}, $\psi_\lambda^k$ is the
$N_J$th term in the power series expansion of \eqref{limpsi2} at $t=0$,
which yields the desired result.
\end{proof}
Define the strictly positive constant $C_J^k$ by
\[C_J^k=\frac{1}{\#W_0}\prod_{\alpha\in\Sigma_{0}^{J,+}}\frac{k_\alpha}{\delta_J(\alpha^\vee)}.
\]
The highest and next to highest homogeneous terms of $p^\lambda_u\in E_J(0)$
($u\in W_0^J$) can now be explicitly computed as follows.
\begin{prop}\label{explicitterms} Let $\lambda\in V_{\mathbb{C}}^*$ be $J$-standard and $u\in W_0^J$.\\
{\bf (i)} The highest homogeneous term $h^\lambda_u$ of
$p^\lambda_u\in E_J(0)$ is of degree $N_J$ and is explicitly given by
\[
h^\lambda_u=C_J^kc_k(u\lambda)\prod_{\alpha\in\Sigma_0^{J,+}}\alpha.
\]
{\bf (ii)} Suppose that $\lambda$ is singular \textup{(}i.e. $J\not=\emptyset$\textup{)}.
The next to highest homogeneous term $n^\lambda_u$ of $p^\lambda_u\in E_J(0)$ is
\[
n^\lambda_u=\partial_{u^{-1}\rho_{u\lambda}^k}\bigl(h^\lambda_u\bigr)
=C_J^kc_k(u\lambda)\sum_{\beta\in\Sigma_0^{J,+}}u\beta(\rho_{u\lambda}^k)
\prod_{\alpha\in\Sigma_0^{J,+}\setminus \{\beta\}}\alpha
\]
with
\begin{equation}\label{rhoexpression}
\rho_\mu^k=\sum_{\alpha\in\Sigma_0^+}\frac{\alpha^\vee}{\mu(\alpha^\vee)+k_\alpha}\in V_{\mathbb{C}}.
\end{equation}
\end{prop}
\begin{rema}
The formula for
$n^\lambda_u$ should be read as an identity between
analytic functions in $k_\alpha>0$
(the possible singularities are easily seen to be removable).
\end{rema}
\begin{proof}
{\bf (i)} Observe that $e_{uv}(0)=e_u(0)$ is independent of $v\in W_{0,J}$, and
\[d_u(0)e_u(0)=c_k(u\lambda)\prod_{\alpha\in\Sigma_{0}^{J,+}}k_\alpha.
\]
Combined with Lemma \ref{differformula} we conclude that the highest homogeneous
term $h^\lambda_u$ of $p^\lambda_u$ is given by
\begin{equation}\label{ht}
\begin{split}
h^\lambda_u&=\frac{C_J^k}{N_J!}c_k(u\lambda)\frac{d^{N_J}}{dt^{N_J}}\bigg|_{t=0}\sum_{v\in W_{0,J}}(-1)^{l(v)}e^{tv\delta_J}\\
&=\frac{C_J^k}{N_J!}c_k(u\lambda)\sum_{v\in
W_{0,J}}(-1)^{l(v)}\bigl(v\delta_J\bigr)^{N_J}.
\end{split}
\end{equation}
On the other hand, by the Weyl denominator formula for $\Sigma_0^J$ we
have
\[
\frac{d^{N_J}}{dt^{N_J}}\bigg|_{t=0}\sum_{v\in W_{0,J}}(-1)^{l(v)}e^{tv\delta_J}=
\frac{d^{N_J}}{dt^{N_J}}\bigg|_{t=0}e^{t\delta_J}\prod_{\alpha\in\Sigma_0^{J,+}}\bigl(1-e^{-t\alpha}\bigr)
=N_J!\prod_{\alpha\in\Sigma_0^{J,+}}\alpha.
\]
Combined with the first equality in \eqref{ht} we obtain the desired expression for $h^\lambda_u$.\\
{\bf (ii)} The next to highest homogeneous term $n^\lambda_u$ of $p^\lambda_u$ is
\[
n^\lambda_u=\frac{N_J}{K_J}\left\{d_u^\prime(0)e_u(0)\sum_{v\in
W_{0,J}}(-1)^{l(v)}\bigl(v\delta_J\bigr)^{N_J-1}
+d_u(0)\sum_{v\in
W_{0,J}}(-1)^{l(v)}e_{uv}^\prime(0)\bigl(v\delta_J\bigr)^{N_J-1}\right\}
\]
in view of Lemma \ref{differformula},
where the prime denotes the $t$-derivative. The first $W_{0,J}$-sum in this expression
is identically zero since it is a $W_{0,J}$-alternating polynomial of degree $<N_J$.
By a direct calculation the remaining expression can be rewritten as
\[n^\lambda_u=\frac{C_J^k}{(N_J-1)!}c_k(u\lambda)\sum_{v\in
W_{0,J}}(-1)^{l(v)}(v\delta_J)(u^{-1}\rho_{u\lambda}^k)\bigl(v\delta_J\bigr)^{N_J-1}.
\]
The desired expression for $n^\lambda_u$ now follows from \eqref{ht}.
\end{proof}
\section{The Bethe ansatz equations}\label{sect8}

In this section we show that $E(\lambda)_Q^{W}\not=\{0\}$
implies that the spectral parameter $\lambda$ is a purely imaginary solution
of the Bethe ansatz equations \eqref{BAE}.

{}From the results of the previous section it is clear that $E(\lambda)_Q^{W}$
is one-dimensional or zero-dimensional. In fact it is one-dimensional if and only if
$Q_k(a_0)\psi_\lambda^k=\psi_\lambda^k$, in which case we have
\[E(\lambda)_Q^{W}=E(\lambda)_Q^{W_0}=\hbox{span}_{\mathbb{C}}\{\psi_\lambda^k\}.
\]
It is convenient to reformulate these observations in terms of
\begin{equation}\label{intertwiner}
\mathcal{J}_k=\partial_{\varphi^\vee}Q_k(a_0)+k_\varphi
\end{equation}
(viewed as an operator on e.g. $C^\infty(V)$ or $E(\lambda)$),
which satisfies the elementary commutation relations
\[\mathcal{J}_k\partial_v=\partial_{s_\varphi v}\mathcal{J}_k,\qquad \forall\,v\in V
\]
(the operator $\mathcal{J}_k$ can be defined on the level of the algebra $H_k$ as the
element $\varphi^\vee\cdot s_0+k_\varphi\in H_k$, in which case it
is the analog of the affine intertwiner from \cite{C1} and \cite[Sect. 4]{O}).
The equality $Q_k(a_0)\psi_\lambda^k=\psi_\lambda^k$ clearly implies $\mathcal{J}_k\psi_\lambda^k=
(\partial_{\varphi^\vee}+k_\varphi)\psi_\lambda^k$.
\begin{lem}\label{regularlem}
If $\lambda$ is regular, then $\mathcal{J}_k\psi_\lambda^k=
(\partial_{\varphi^\vee}+k_\varphi)\psi_\lambda^k$ implies $Q_k(a_0)\psi_\lambda^k=\psi_\lambda^k$.
\end{lem}
\begin{proof}
By \eqref{Steinbergregular} we have a unique expansion
\[Q_k(a_0)\psi_\lambda^k-\psi_\lambda^k=\sum_{w\in W_0}d_we^{w\lambda}
\]
with $d_w\in\mathbb{C}$. We conclude from the equality
$\mathcal{J}_k\psi_\lambda^k=(\partial_{\varphi^\vee}+k_\varphi)\psi_\lambda^k$ that
$w\lambda(\varphi^\vee)d_w=0$ for all $w\in W_0$. Since $\lambda$ is regular, this implies $d_w=0$ for all $w\in W_0$.
\end{proof}

For $p\in P(V)_{\mathbb{C}}\simeq S(V^*)_{\mathbb{C}}$ we write $p(\partial^\mu)$ for
the associated constant coefficient differential
operator acting on smooth functions in $\mu\in V_{\mathbb{C}}^*$.
\begin{lem}\label{techlem}
Let $p\in P(V)_{\mathbb{C}}\simeq S(V^*)_{\mathbb{C}}$.
For $w\in W_0$ we have
\begin{equation*}
\begin{split}
\mathcal{J}_k\bigl(p(w^{-1}\cdot)e^{w\mu}\bigr)(v)&=
-p(\partial^\mu)\bigl((\mu(w^{-1}\varphi^\vee)+k_\varphi)e^{\mu(w^{-1}\varphi^\vee)}e^{\mu(w^{-1}s_{\varphi}v)}\bigr),\\
\bigl(\partial_{\varphi^\vee}+k_{\varphi}\bigr)\bigl(p(w^{-1}\cdot)e^{w\mu}\bigr)(v)&=
p(\partial^\mu)\bigl((\mu(w^{-1}\varphi^\vee)+k_\varphi)e^{\mu(w^{-1}v)}\bigr),
\end{split}
\end{equation*}
where we view the left hand sides as functions in $v\in V$
and the right hand sides as functions in $\mu\in V_{\mathbb{C}}^*$.
In particular,
\[\mathcal{J}_k\bigl(P^{(\leq j)}(V)_{\mathbb{C}}\,e^\mu\bigr)\subseteq
P^{(\leq j)}(V)_{\mathbb{C}}\,e^{s_\varphi\mu},\qquad
\bigl(\partial_{\varphi^\vee}+k_\varphi\bigr)\bigl(P^{(\leq j)}(V)_{\mathbb{C}}\,
e^\mu\bigr)\subseteq P^{(\leq j)}(V)_{\mathbb{C}}\,e^{\mu}
\]
for $j\in\mathbb{Z}_{\geq 0}$ and $\mu\in V_{\mathbb{C}}^*$.
\end{lem}
\begin{proof}
Observe that
\[\bigl(p(w^{-1}\cdot)e^{w\mu}\bigr)(v)=p(\partial^\mu)\bigl(e^{\mu(w^{-1}v)}\bigr),
\]
and $p(\partial^\mu)$ (acting on $\mu\in V_{\mathbb{C}}^*$) clearly commutes with $\mathcal{J}_k$
and $(\partial_{\varphi^\vee}+k_\varphi)$ (which act on $v\in V$).
Thus it suffices to prove the lemma for $p\equiv 1$,
in which case the second formula is trivial.
To prove the first formula for $p\equiv 1$ we may assume
without loss of generality that $w=e$ is the unit element of $W_0$.
Suppose that $\mu\in V_{\mathbb{C}}^*$ is regular. A direct computation
using the definition \eqref{Q}
of $Q_k(a_0)$ as an integral-reflection operator yields
\[Q_k(a_0)e^{\mu}=-\frac{k_\varphi}{\mu(\varphi^\vee)}e^\mu+\left(\frac{\mu(\varphi^\vee)+k_\varphi}{\mu(\varphi^\vee)}\right)
e^{\mu(\varphi^\vee)}e^{s_\varphi\mu},
\]
hence
\[\mathcal{J}_k(e^\mu)=-\bigl(\mu(\varphi^\vee)+k_\varphi\bigr)e^{\mu(\varphi^\vee)}e^{s_\varphi\mu}.
\]
In the latter formula the regularity constraint on $\mu$ can be removed by continuity.
\end{proof}

We denote $\pi_\lambda^{(j)}: E(\lambda)\rightarrow
E^{(j)}(\lambda)$ for the projection onto $E^{(j)}(\lambda)$
along the decomposition $E(\lambda)=\bigoplus_{r=0}^\infty
E^{(r)}(\lambda)$. Observe that
\begin{equation}\label{identitydecomposition}
\hbox{Id}_{E(\lambda)}=\sum_{j=0}^{N_J}\pi_\lambda^{(j)}
\end{equation}
if $\lambda$ is $J$-standard in view of
Proposition \ref{explicitterms}{\bf (i)}.
In this section we consider the constraint on $\lambda$ such that
\begin{equation}\label{jrestraint}
\pi_\lambda^{(j)}\bigl(\mathcal{J}_k\psi_\lambda^k\bigr)=
\pi_\lambda^{(j)}\bigl((\partial_{\varphi^\vee}+k_\varphi)\psi_\lambda^k\bigr)
\end{equation}
for the highest degree component $j=N_J$.

The map $u\mapsto u^J$, where $u^J\in W_0^J$ is obtained from
the unique decomposition
\begin{equation}\label{invo}
s_{\varphi}u=u^Ju_J,\qquad u^J\in W_0^J,\,\, u_J\in W_{0,J},
\end{equation}
defines an involution on $W_0^J$. Observe that
\begin{equation}\label{invo2}
(u^J)_J=(u_J)^{-1},\qquad u\in W_0^J.
\end{equation}
Recall that $c_k$ denotes the regularized $c$-function \eqref{regularc}.
\begin{lem}\label{highestBAElem}
Suppose that $\lambda\in V_{\mathbb{C}}^*$ is $J$-standard.\\
{\bf (i)}
The equation \eqref{jrestraint} for $j=N_J$ holds if and only if
$\lambda$ satisfies the equations
\begin{equation}\label{parts}
c_k(s_{\varphi}u\lambda)(u\lambda(\varphi^\vee)-k_\varphi)e^{-u\lambda(\varphi^\vee)}(-1)^{l(u_J)}=
c_k(u\lambda)(u\lambda(\varphi^\vee)+k_\varphi),\quad \forall\,u\in
W_0^J.
\end{equation}
{\bf (ii)} For $u\in W_{0}^J$ and for multiplicity functions $k$ such that $c_k(u\lambda)\not=0$, we have
\[\frac{c_k(s_{\varphi}u\lambda)}{c_k(u\lambda)}=(-1)^{l(u_J)}\prod_{\alpha\in\Sigma_0^+\cap s_\varphi\Sigma_0^-}
\frac{u\lambda(\alpha^\vee)-k_\alpha}{u\lambda(\alpha^\vee)+k_\alpha}.
\]
\end{lem}
\begin{proof}
{\bf (i)} By \eqref{exp}, Lemma \ref{techlem} and Proposition \ref{explicitterms}{\bf (i)} we have
\begin{equation}\label{tt1}
\begin{split}
\pi_\lambda^{(N_J)}(\mathcal{J}_k\psi_\lambda^k)&=-C_J^k\sum_{u\in
W_0^J}c_k(u\lambda)(u\lambda(\varphi^\vee)+k_\varphi)e^{u\lambda(\varphi^\vee)}
e^{s_{\varphi}u\lambda}\prod_{\alpha\in\Sigma_0^{J,+}}s_{\varphi}u\alpha,\\
\pi_\lambda^{(N_J)}\bigl((\partial_{\varphi^\vee}+k_\varphi)\psi_\lambda^k\bigr)&=
C_J^k\sum_{u\in
W_0^J}c_k(u\lambda)(u\lambda(\varphi^\vee)+k_\varphi)e^{u\lambda}\prod_{\alpha\in\Sigma_0^{J,+}}u\alpha.
\end{split}
\end{equation}
The proof now follows by equating the coefficients of $e^{u\lambda}\prod_{\alpha\in\Sigma_0^{J,+}}u\alpha$
($u\in W_0^J$) in \eqref{tt1} using \eqref{invo}.\\
{\bf (ii)} We first compare the denominators of $c_k(u\lambda)$ and $c_k(s_\varphi u\lambda)=c_k(u^J\lambda)$.
If $\mu\in V_{\mathbb{C}}^*$ is regular then
\begin{equation*}
\begin{split}
\prod_{\alpha\in\Sigma_0^+\setminus u^J\Sigma_0^{J,+}}u^J\mu(\alpha^\vee)&=
\prod_{\alpha\in\Sigma_0^+}u^J\mu(\alpha^\vee)\prod_{\beta\in uu_J^{-1}\Sigma_0^{J,+}}(uu_J^{-1}\mu(\beta^\vee))^{-1}\\
&=(-1)^{l(u_J)}\prod_{\alpha\in\Sigma_0^+}s_\varphi uu_J^{-1}\mu(\alpha^\vee)
\prod_{\beta\in u\Sigma_0^{J,+}}(uu_J^{-1}\mu(\beta^\vee))^{-1}\\
&=
(-1)^{l(u_J)+1}\prod_{\alpha\in\Sigma_0^+\setminus u\Sigma_0^{J,+}}uu_J^{-1}\mu(\alpha^\vee).
\end{split}
\end{equation*}
Taking the limit $\mu\rightarrow \lambda$ we obtain
\[\prod_{\alpha\in\Sigma_0^+\setminus u^J\Sigma_0^{J,+}}u^J\lambda(\alpha^\vee)=(-1)^{l(u_J)+1}\prod_{\alpha\in\Sigma_0^+\setminus
u\Sigma_0^{J,+}}u\lambda(\alpha^\vee).
\]
A similar (and easier) computation leads to the comparative formula
\[\prod_{\alpha\in \Sigma_0^+\setminus u^J\Sigma_0^{J,+}}\bigl(u^J\lambda(\alpha^\vee)+k_\alpha\bigr)=
-\left(\prod_{\beta\in\Sigma_0^+\cap s_\varphi\Sigma_0^-}\frac{u\lambda(\beta^\vee)-k_\beta}{u\lambda(\beta^\vee)+k_\beta}\right)
\prod_{\alpha\in \Sigma_0^+\setminus u\Sigma_0^{J,+}}\bigl(u\lambda(\alpha^\vee)+k_\alpha\bigr)
\]
for the numerators of $c_k(u\lambda)$ and $c_k(u^J\lambda)$. Combining both formulas leads to the desired result.
\end{proof}
Recall from Section \ref{sect2} that $\hbox{BAE}_k$ is the set of purely imaginary solutions of the Bethe ansatz equations
\eqref{BAE}.
\begin{prop}\label{highestBAE}
Suppose that $\lambda\in V_{\mathbb{C}}^*$ is $J$-standard.
The equation \eqref{jrestraint} for $j=N_J$ holds if and only if
$\lambda\in \textup{BAE}_k$.
\end{prop}
\begin{proof}
We first show that $\lambda$ is purely imaginary if $\lambda$ satisfies the equation \eqref{parts}.
Let $\mu=u\lambda$ ($u\in W_0^J$) be the element in the $W_0$-orbit of $\lambda$ having its real part in $\overline{V_+^*}$.
Then $c_k(\mu)\not=0$ since the multiplicity function $k$ is strictly positive,
hence \eqref{parts} and Lemma \ref{highestBAElem}{\bf (ii)} imply
\begin{equation}\label{im}
e^{\mu(\varphi^\vee)}=\frac{\mu(\varphi^\vee)-k_\varphi}{\mu(\varphi^\vee)+k_\varphi}
\prod_{\alpha\in\Sigma_0^+\cap
s_{\varphi}\Sigma_0^-}\frac{\mu(\alpha^\vee)-
k_\alpha}{\mu(\alpha^\vee)+k_\alpha}.
\end{equation}
The modulus of the left hand
(respectively right hand side) of \eqref{im} is $\geq 1$
(respectively $\leq 1$) since the real part of $\mu$ is in $\overline{V_+^*}$ and
the multiplicity function $k$ is strictly positive. Thus $|e^{\mu(\varphi^\vee)}|=1$, implying that
$\mu(\varphi^\vee)$ is purely imaginary. Since $\varphi^\vee=\sum_{j=1}^nm_ja_j^\vee$
with $m_j$ strictly positive integers and since the real part of $\mu$ lies in $\overline{V_+^*}$,
we conclude that $\mu(a_j^\vee)$ is purely imaginary
for all co-roots $a_j^\vee$ ($j=1,\ldots,n$). This implies $\mu\in iV^*$, hence $\lambda\in iV^*$.

Combined with Lemma \ref{highestBAElem}{\bf (i)} it follows that $\lambda$ satisfies \eqref{jrestraint} for $j=N_J$
if and only if $\lambda$ is a purely imaginary solution of the equations \eqref{parts}. For purely imaginary $\lambda$
we have $c_k(u\lambda)\not=0$ for all $u\in W_0^J$ due to the strict positivity of the multiplicity function $k$. The proof
now follows from Lemma \ref{highestBAElem}{\bf (ii)} and Remark \ref{BAEremark}.
\end{proof}
As an immediate result we obtain the following ``regular part'' of Theorem \ref{main2}.
\begin{cor}\label{regularcor}
Suppose that $\lambda\in V_{\mathbb{C}}^*$ is regular. The space $E(\lambda)_Q^{W}$ is zero-dimensional or one-dimensional.
It is one-dimensional if and only if $\lambda\in \textup{BAE}_k$.
In that case  $E(\lambda)_Q^{W}$ is spanned by $\psi_\lambda^k$ \eqref{BAF2}.
\end{cor}
\begin{proof}
By the observations at the beginning of the section it suffices to show that $E(\lambda)_Q^W\not=\{0\}$ iff
$\lambda\in\hbox{BAE}_k$.

Since $\hbox{BAE}_k\subset iV^*$ is a $W_0$-invariant subset and $E(\lambda)_Q^{W}$ only depends on the $W_0$-orbit
of $\lambda$, we may assume without loss of generality that $\lambda$ is $\emptyset$-standard.
If $E(\lambda)_Q^{W}\not=\{0\}$ then \eqref{jrestraint} holds, hence $\lambda\in \hbox{BAE}_k$ by Proposition \ref{highestBAE}.
Conversely, suppose that $\lambda\in \hbox{BAE}_k$. Since $\lambda$ is regular we
have $\hbox{Id}_{E(\lambda)}=\pi_\lambda^{(0)}$ by \eqref{identitydecomposition}, hence
$\mathcal{J}_k\psi_\lambda^k=\bigl(\partial_{\varphi^\vee}+k_\varphi\bigr)\psi_\lambda^k$ by Proposition \ref{highestBAE}.
By Lemma \ref{regularlem} this implies $Q_k(a_0)\psi_\lambda^k=\psi_\lambda^k$, hence $0\not=\psi_\lambda^k\in E(\lambda)_Q^{W}$.
\end{proof}

\section{The master function}\label{sect9}
In this section we prove Proposition \ref{parametrization}, which yields a parametrization
of the set $\hbox{BAE}_k$ of purely imaginary solutions of the Bethe ansatz equations \eqref{BAE} by the weight lattice $P$.

We first rewrite the Bethe ansatz equations \eqref{BAE} in logarithmic form.
By a direct computation using the elementary identity
\[e^{-2i\arctan(x)}=\frac{1-ix}{1+ix}\qquad (x\in \mathbb{R})
\]
the Bethe ansatz equations \eqref{BAE} for $\lambda\in iV^*$ can be rewritten as
\begin{equation}\label{for1}
-i\lambda(w\varphi^\vee)+\sum_{\alpha\in\Sigma_0}
\arctan\left(\frac{-i\lambda(\alpha^\vee)}{k_\alpha}\right)\alpha(w\varphi^\vee)=0
\quad \hbox{ modulo }\,2\pi\mathbb{Z}
\end{equation}
for all $w\in W_0$. On the other hand, for $\mu\in P$ the gradient of the master function
$S_k(\mu,\cdot): V^*\rightarrow\mathbb{R}$ (see \eqref{master}) is determined by
\begin{equation}\label{for2}
\bigl(\partial_{\xi}S_k(\mu,\cdot)\bigr)(\eta)=\langle
\eta-2\pi\mu+\sum_{\alpha\in\Sigma_0}\arctan\left(\frac{\eta(\alpha^\vee)}{k_\alpha}\right)\alpha,\xi\rangle, \qquad \xi, \eta\in V^*
\end{equation}
Comparing \eqref{for1} and \eqref{for2} yields the following result.
\begin{lem}\label{steppar}
We have $\lambda\in \textup{BAE}_k$ if and only if $\lambda=i\eta$ with $\eta\in V^*$
an extremal vector of the master function $S_k(\mu,\cdot)$ for some $\mu\in P$.
\end{lem}
\begin{proof}
$\Sigma_0$ is an irreducible root system in $V^*$, hence $\{w\varphi \, | \, w\in W_0\}$
spans $V^*$. Thus $\eta\in V^*$ is an extremal
vector of $S_k(\mu,\cdot)$ if and only if $\bigl(\partial_{w\varphi}S_k(\mu,\cdot)\bigr)(\eta)=0$ for all $w\in W_0$,
which by \eqref{for2} is equivalent to
\[\eta(w\varphi^\vee)+\sum_{\alpha\in\Sigma_0}\alpha(w\varphi^\vee)\arctan\left(\frac{\eta(\alpha^\vee)}{k_\alpha}\right)=2\pi\mu(w\varphi^\vee)
\]
for all $w\in W_0$. Comparing to \eqref{for1},
the proof now follows from \eqref{Palt}.
\end{proof}
We thus need to analyze the extrema of the master function $S_k(\mu,\cdot)$ at a given weight $\mu\in P$.
Observe that the Hessian $B_\xi^k: V^*\times V^*\rightarrow \mathbb{R}$
of $S_k(\mu,\cdot)$ at $\xi\in V^*$ is independent of $\mu$,
and is given explicitly by
\begin{equation}\label{Hessian}
\begin{split}
B_\xi^k(\eta,\eta^\prime)&=\left(\partial_\eta\partial_{\eta^{\prime}}S_k(\mu,\cdot)
\right)(\xi)\\
&=\langle\eta,\eta^\prime\rangle+\frac{1}{2}\sum_{\alpha\in\Sigma_0}
k_\alpha\|\alpha\|^2\frac{\eta(\alpha^\vee)\eta^\prime(\alpha^\vee)}
{k_\alpha^2+\xi(\alpha^\vee)^2},\qquad \eta,\eta^\prime\in V^*.
\end{split}
\end{equation}
By the strict positivity of the multiplicity function $k$, it follows from \eqref{Hessian} that
the Hessian $B_\xi^k$ is positive definite for all $\xi\in V^*$, hence $S_k(\mu,\cdot)$ is strictly convex. Furthermore, for all $\mu\in P$,
\[S_k(\mu,\xi)\geq \frac{\|\xi\|^2}{2}-2\pi\langle \mu,\xi\rangle\rightarrow \infty,\qquad\|\xi\|\rightarrow\infty
\]
hence $S_k(\mu,\cdot)$ has a unique extremum $\widehat{\mu}_k\in V^*$,
which is a global minimum. It now follows from \eqref{for2} that $\widehat{\mu}_k$ ($\mu\in P$) is uniquely determined by the equation
\begin{equation}
\widehat{\mu}_k+\sigma_{\widehat{\mu}_k}^k=2\pi\mu
\end{equation}
in $V^*$, where $\sigma_\lambda^k\in V^*$ ($\lambda\in V^*$) is defined by
\[\sigma_\lambda^k=\sum_{\alpha\in\Sigma_0}
\arctan\left(\frac{\lambda(\alpha^\vee)}{k_\alpha}\right)\alpha.
\]
Combined with Lemma \ref{steppar} it now follows that the
map $\mu\mapsto i\widehat{\mu}_k$ is a bijection from the weight lattice $P$ onto $\hbox{BAE}_k$. The $W_0$-equivariance of
this map is immediate from the equivariance property
\[\bigl(\partial_{w\xi}S_k(w\mu,\cdot)\bigr)(w\eta)=\bigl(\partial_{\xi}S_k(\mu,\cdot)\bigr)(\eta),\qquad \forall\,w\in W_0
\]
for $\xi,\eta\in V^*$ and $\mu\in P$. This completes the proof of Proposition \ref{parametrization}.


\section{Moment gaps}\label{sect10}
In this section we prove Proposition \ref{mg}, which yields estimates for
the location of the deformed weight $\widehat{\mu}=\widehat{\mu}_k$ compared to the
parametrizing weight $\mu\in P$. In view of \eqref{for2} and Lemma \ref{steppar},
the deformed weight $\widehat{\mu}\in V^*$ ($\mu\in P$) is the unique solution of

The following lemma establishes the necessary bounds for $\sigma_\lambda^k$.
\begin{lem}\label{rholem}
For $\lambda\in \overline{V_+^*}$,
\[0\leq \sigma_\lambda^k(\beta^\vee)\leq
\frac{h_k}{n}\lambda(\beta^\vee),\qquad \forall\,\beta\in\Sigma_0^+
\]
with $h_k=2\sum_{\alpha\in\Sigma_0}k_\alpha^{-1}$.
\end{lem}
\begin{proof}
Fix $\lambda\in \overline{V_+^*}$ and $\beta\in\Sigma_0^+$.
Let $\Sigma_0^\beta$ be the set of roots $\alpha\in\Sigma_0$ satisfying $\alpha(\beta^\vee)>0$, then
\begin{equation}\label{convenformrho}
\sigma_\lambda^k(\beta^\vee)=\sum_{\alpha\in\Sigma_0^\beta}\left\{
\arctan\left(\frac{\lambda(\alpha^\vee)}{k_\alpha}\right)-
\arctan\left(\frac{\lambda(s_\beta\alpha^\vee)}{k_\alpha}\right)
\right\}\alpha(\beta^\vee).
\end{equation}
Each term in this sum is positive, hence $\sigma_\lambda^k(\beta^\vee)\geq 0$.

For the second inequality, we use the estimate for $\alpha\in\Sigma_0^\beta$,
\[\arctan\left(\frac{\lambda(\alpha^\vee)}{k_\alpha}\right)-
\arctan\left(\frac{\lambda(s_\beta(\alpha^\vee))}{k_\alpha}\right)=
\int_{\lambda(s_\beta(\alpha^\vee))/k_\alpha}^{\lambda(\alpha^\vee)/k_\alpha}\frac{dx}{1+x^2}
\leq \frac{\lambda(\beta^\vee)\beta(\alpha^\vee)}{k_\alpha},
\]
leading to
\begin{equation}\label{almostineq}
\sigma_\lambda^k(\beta^\vee)\leq \lambda(\beta^\vee)
\sum_{\alpha\in\Sigma_0^\beta}\frac{\beta(\alpha^\vee)\alpha(\beta^\vee)}{k_\alpha}
=\frac{\lambda(\beta^\vee)}{2}\sum_{\alpha\in\Sigma_0}\frac{\beta(\alpha^\vee)\alpha(\beta^\vee)}{k_\alpha}
\end{equation}
in view of \eqref{convenformrho}. Now note that
\[\xi\mapsto\sum_{\alpha\in\Sigma_0}k_\alpha^{-1}\xi(\alpha^\vee)\alpha
\]
defines a $W_0$-equivariant linear map $V^*\rightarrow V^*$. By
Schur's lemma it equals $C_k\hbox{Id}_{V^*}$ for some constant
$C_k\in\mathbb{C}$. To determine $C_k$ explicitly we fix a basis $\{e_j\}_{j=1}^n$
of $V$ and we denote $\{\epsilon_j\}_{j=1}^n$ for the corresponding
dual basis of $V^*$. Then
\[C_k n=\sum_{j=1}^n\sum_{\alpha\in\Sigma_0}k_\alpha^{-1}\epsilon_j(\alpha^\vee)\alpha(e_j)
=h_k
\]
with $h_k=2\sum_{\alpha\in\Sigma_0}k_\alpha^{-1}$. Combined with
\eqref{almostineq} we obtain $\sigma_\lambda^k(\beta^\vee)\leq \frac{h_k}{n}\lambda(\beta^\vee)$.
\end{proof}
\begin{cor}\label{rhocor}
Let $\mu\in P$. We have $\widehat{\mu}_k\in \overline{V_+^*}$ if and only if $\mu\in P^+$.
\end{cor}
\begin{proof}
Let $\mu\in P$ and suppose that $\widehat{\mu}_k\in
\overline{V_+^*}$. Then for all $\beta\in\Sigma_0^+$,
\[2\pi\mu(\beta^\vee)=\widehat{\mu}_k(\beta^\vee)+\sigma_{\widehat{\mu}_k}^k(\beta^\vee)\geq
0\]
by Lemma \ref{rholem}, hence $\mu\in P^+$.

Conversely, suppose that $\mu\in P^+$ and let $w\in W_0$ such that
$w\widehat{\mu}_k\in \overline{V_+^*}$. By Proposition
\ref{parametrization} this implies $\widehat{w\mu}_k\in
\overline{V_+^*}$. By the previous paragraph we conclude that $w\mu\in
P^+$. On the other hand $P^+\cap W_0\mu=\{\mu\}$, hence $w\mu=\mu\in P^+$ and
$\widehat{\mu}_k=\widehat{w\mu}_k\in\overline{V_+^*}$.
\end{proof}
Proposition \ref{mg} is now a direct consequence of Corollary \ref{rhocor}
and Lemma \ref{rholem}.


\section{The Pauli principle}\label{sect11}

In this section we complete the proof of Theorem \ref{main2} (and hence also of Theorem \ref{main}).
In view of Proposition \ref{highestBAE} and Corollary \ref{regularcor}
it suffices to show the following root system analog of the Pauli principle.
\begin{prop}\label{Pauli}
If $\lambda\in \textup{BAE}_k$ is singular then $E(\lambda)_Q^{W}=\{0\}$.
\end{prop}
For the proof of Proposition \ref{Pauli} we may assume without loss of generality that
$\lambda\in\textup{BAE}_k$ is $J$-standard (in particular, $\lambda\in i\overline{V_+^*}$).
We write $V_J^*\subseteq V^*$ for the real sub-space
spanned by the subset $J$ of simple roots. Its complement in $V$
is defined by
\[V_J^\perp=\{v\in V \, | \, \xi(v)=0\quad \forall\,\xi\in V_J^*\}.
\]
Observe that $V_J^{\perp}=V$ iff $J=\emptyset$ iff $\lambda$ is regular.

Consider the linear map $K_\lambda^k: V\rightarrow V$ defined by
\[K_{\lambda}^k(v)=v+\sum_{\alpha\in\Sigma_0}
\frac{k_\alpha\alpha(v)\alpha^\vee}{k_\alpha^2-\lambda(\alpha^\vee)^2},
\qquad v\in V.
\]

\begin{lem}\label{Kform}
Let $\lambda\in iV^*$ be a singular $J$-standard solution of the
Bethe ansatz equations \eqref{BAE}. Then $\lambda$ satisfies the constraint
\begin{equation}\label{restraintonelower}
\pi_\lambda^{(N_J-1)}\bigl(\mathcal{J}_k\psi_\lambda^k\bigr)=
\pi_\lambda^{(N_J-1)}\bigl((\partial_{\varphi^\vee}+
k_\varphi)\psi_\lambda^k\bigr)
\end{equation}
iff $K_\lambda^k(V)\subseteq V_J^{\perp}$.
\end{lem}
\begin{proof}
Fix a singular $J$-standard solution $\lambda\in i\overline{V_+^*}$ of the Bethe ansatz equations \eqref{BAE}
(in particular $J\not=\emptyset$). By a similar computation as in the proof of Proposition \ref{highestBAE}
we obtain from \eqref{exp}, Lemma \ref{techlem} and
Proposition \ref{explicitterms},
\begin{equation*}
\begin{split}
\pi_\lambda^{(N_J-1)}\bigl((\partial_{\varphi^\vee}+k_\varphi)\psi_\lambda^k\bigr)&=
C_J^k\sum_{u\in
W_0^J}c_k(u\lambda)\sum_{\beta\in\Sigma_0^{J,+}}u\beta(a_{u\lambda})e^{u\lambda}\prod_{\alpha\in\Sigma_0^{J,+}\setminus
\{\beta\}}u\alpha,\\
\pi_\lambda^{(N_J-1)}\bigl(\mathcal{J}_k\psi_\lambda\bigr)&=
C_J^k\sum_{u\in W_0^J}c_k(u\lambda)e^{-u^J\lambda(\varphi^\vee)}
\sum_{\beta\in\Sigma_0^{J,+}}u\beta(b_{u^J\lambda})e^{u^J\lambda}\prod_{\alpha\in\Sigma_0^{J,+}\setminus
\{\beta\}}u^Ju_J\alpha
\end{split}
\end{equation*}
with vectors $a_\mu,b_\mu\in V_{\mathbb{C}}$ ($\mu\in V_{\mathbb{C}}^*$) given by
\begin{equation*}
\begin{split}
a_\mu&=(\mu(\varphi^\vee)+k_\varphi)\rho_{\mu}^k+\varphi^\vee,\\
b_\mu&=\bigl(\mu(\varphi^\vee)-k_\varphi\bigr)
\bigl(\rho_{s_{\varphi}\mu}^k+\varphi^\vee\bigr)-\varphi^\vee,
\end{split}
\end{equation*}
where we have used the involution on $W_0^J$ defined by \eqref{invo}, as well as \eqref{invo2}.
For $u\in W_0^J$ we have
\begin{equation*}
\begin{split}
\sum_{\beta\in\Sigma_0^{J,+}}u\beta(b_{u^J\lambda})\prod_{\alpha\in\Sigma_0^{J,+}\setminus
\{\beta\}}u^Ju_J\alpha&=(-1)^{l(u_J)}\left(\sum_{\beta\in\Sigma_0^{J,+}}
\frac{u\beta(b_{u^J\lambda})}{u^Ju_J\beta}\right)\prod_{\alpha\in\Sigma_0^{J,+}}u^J\alpha\\
&=\frac{1}{2}(-1)^{l(u_J)}\left(\sum_{\beta\in\Sigma_0^J}\frac{uu_J^{-1}\beta(b_{u^J\lambda})}{u^J\beta}\right)
\prod_{\alpha\in\Sigma_0^{J,+}}u^J\alpha\\
&=(-1)^{l(u_J)}\sum_{\beta\in\Sigma_0^{J,+}}uu_J^{-1}\beta(b_{u^J\lambda})
\prod_{\alpha\in\Sigma_0^{J,+}\setminus \{\beta\}}u^J\alpha.
\end{split}
\end{equation*}
Consequently \eqref{restraintonelower} is equivalent to
\[ c_k(u\lambda)u\beta(a_{u\lambda})=(-1)^{l(u_J)}c_k(u^J\lambda)e^{-u\lambda(\varphi^\vee)}s_{\varphi}u\beta(b_{u\lambda}),\qquad
\forall\, u\in W_0^J,\,\,\forall\, \beta\in\Sigma_0^{J,+}.
\]
Since $\lambda$ is a solution of the Bethe ansatz equations (see \eqref{parts} for the convenient equivalent form
of the Bethe ansatz equations) this is equivalent to
\begin{equation}\label{formulas}
\bigl(u\lambda(\varphi^\vee)-k_\varphi\bigr)a_{u\lambda}-
\bigl(u\lambda(\varphi^\vee)+k_\varphi\bigr)s_{\varphi}b_{u\lambda}\in u(V_J^\perp),
\qquad \forall\, u\in W_0^J.
\end{equation}
Note that \eqref{formulas} only depends on the coset $uW_{0,J}$ ($u\in
W_0^J$). Using the explicit expressions for $a_{u\lambda}$ and $b_{u\lambda}$ we can rewrite \eqref{formulas}
as
\begin{equation}\label{formulasone}
\bigl(w^{-1}\rho_{w\lambda}^k-w^{-1}s_{\varphi}\rho_{s_{\varphi}w\lambda}^k\bigr)
+\left(\frac{w\lambda(\varphi^\vee)^2-k_{\varphi}^2-2k_{\varphi}}{w\lambda(\varphi)^2-k_{\varphi}^2}\right)
w^{-1}\varphi^\vee\in V_J^\perp,\quad \forall\,w\in W_0.
\end{equation}
We match \eqref{formulasone} to the desired condition $K_\lambda^k(V)\subseteq V_J^\perp$ as follows.
Since $\Sigma_0$ is an irreducible root system in $V^*$, the condition $K_\lambda^k(V)\subseteq V_J^\perp$
is equivalent to $K_\lambda^k(w^{-1}\varphi^\vee)\in V_J^\perp$ for all $w\in W_0$, which in turn is
equivalent to \eqref{formulasone} if
\begin{equation}\label{rhoK}
K_\lambda^k(w^{-1}\varphi^\vee)=\bigl(w^{-1}\rho_{w\lambda}^k-w^{-1}s_{\varphi}\rho_{s_{\varphi}w\lambda}^k\bigr)
+\left(\frac{w\lambda(\varphi^\vee)^2-k_{\varphi}^2-2k_{\varphi}}{w\lambda(\varphi)^2-k_{\varphi}^2}\right)
w^{-1}\varphi^\vee
\end{equation}
for all $w\in W_0$. To prove \eqref{rhoK} we first observe that
\[s_\varphi\rho_{s_{\varphi}w\lambda}^k=\rho_{w\lambda}^k-
2\sum_{\alpha\in\Sigma_0^+\cap s_{\varphi}\Sigma_0^-}\frac{k_\alpha\alpha^\vee}{k_\alpha^2-w\lambda(\alpha^\vee)^2}
\]
by the explicit expression \eqref{rhoexpression} for $\rho_\mu^k$. Using \eqref{highestrootproperty} this can be rewritten as
\[
w^{-1}\rho_{w\lambda}^k-w^{-1}s_{\varphi}\rho_{s_{\varphi}w\lambda}^k=
2\frac{k_{\varphi}w^{-1}\varphi^\vee}{w\lambda(\varphi^\vee)^2-k_{\varphi}^2}+
2\sum_{\alpha\in\Sigma_0^+}\frac{k_\alpha\alpha(\varphi^\vee)w^{-1}\alpha^\vee}{k_\alpha^2-w\lambda(\alpha^\vee)^2}.
\]
The second term can be rewritten as
\begin{equation*}
\begin{split}
2\sum_{\alpha\in\Sigma_0^+}\frac{k_\alpha\alpha(\varphi^\vee)w^{-1}\alpha^\vee}{k_\alpha^2-w\lambda(\alpha^\vee)^2}&=
\sum_{\alpha\in\Sigma_0}\frac{k_\alpha\alpha(\varphi^\vee)w^{-1}\alpha^\vee}{k_\alpha^2-w\lambda(\alpha^\vee)^2}\\
&=\sum_{\alpha\in\Sigma_0}\frac{k_\alpha\alpha(w^{-1}\varphi^\vee)}{k_\alpha^2-\lambda(\alpha^\vee)^2}\\
&=K_\lambda^k(w^{-1}\varphi^\vee)-w^{-1}\varphi^\vee.
\end{split}
\end{equation*}
Combining the latter two formulas yields \eqref{rhoK}.
\end{proof}

It follows from \eqref{Hessian} that
\[B_{-i\lambda}^k(\eta_v,\eta_{v^\prime})=\langle K_\lambda^k(v),v^\prime\rangle,\qquad v,v^\prime\in V
\]
with $\eta_v=\langle v,\cdot\rangle\in V^*$ and
$B_{-i\lambda}^k$ the Hessian of the master function $S_k$ at $-i\lambda\in V^*$.
Since $B_{-i\lambda}^k$ is positive definite, $K_\lambda^k: V\overset{\sim}{\longrightarrow} V$ is a linear
isomorphism. Proposition \ref{Pauli} thus is an immediate consequence of Lemma \ref{Kform}.


\end{document}